\documentclass[a4paper,11pt, reqno]{amsart} 
\usepackage{amssymb,amsthm,amsmath}

\usepackage{amsfonts}
\usepackage{amscd}
\usepackage{amssymb}
\usepackage{enumerate}
\usepackage{graphicx}
\allowdisplaybreaks
\usepackage{mathabx}
\usepackage{color}
\usepackage{amsbsy}
\usepackage{graphicx}
\usepackage{amsthm}
\usepackage{amsmath}
\usepackage{amsxtra}
\usepackage{mathrsfs}
\usepackage{bbm}
\usepackage{dsfont}

\usepackage{ifthen}

\usepackage[colorlinks,citecolor=red,pagebackref,hypertexnames=false]{hyperref} 

\newtheorem{theorem}{Theorem}[section]
\newtheorem{corollary}[theorem]{Corollary}

\newtheorem{prop}[theorem]{Proposition}

\theoremstyle{definition}
\newtheorem{defn}[theorem]{Definition}
\theoremstyle{remark}
\newtheorem{rem}[theorem]{Remark}
\numberwithin{equation}{section}
\definecolor{red}{rgb}{1.0, 0.0, 0.0}
\newcommand{\Bea}{\begin{eqnarray*}}
	\newcommand{\Eea}{\end{eqnarray*}}
\newcommand{\Be} {\begin{equation*}}
	\newcommand{\Ee} {\end{equation*}}
\newcommand{\be} {\begin{equation}}
	\newcommand{\ee} {\end{equation}}
\newcommand{\bea} {\begin{eqnarray}}
	\newcommand{\eea} {\end{eqnarray}}

\title[ Multilinear Fourier Integral operators on modulation spaces ]{ Multilinear Fourier Integral operators on modulation spaces }
\author[Aparajita Dasgupta]{Aparajita Dasgupta}
\address{
	Aparajita Dasgupta:
	\endgraf
	Department of Mathematics
	\endgraf
	Indian Institute of Technology, Delhi, Hauz Khas
	\endgraf
	New Delhi-110016 
	\endgraf
	India
	\endgraf
	{\it E-mail address:} {\rm adasgupta@maths.iitd.ac.in}
}
\author[Lalit Mohan]{Lalit Mohan}
\address{
	Lalit Mohan:
	\endgraf
	Department of Mathematics
	\endgraf
	Indian Institute of Technology, Delhi
	\endgraf
	India
	\endgraf
	{\it E-mail address:} {\rm mohanlalit871@gmail.com}
}
\author[Shyam Swarup Mondal]{Shyam Swarup Mondal}
\address{
	Shyam Swarup Mondal:
	\endgraf
	Department of Mathematics
	\endgraf
	Indian Institute of Technology, Delhi
	\endgraf
	India
	\endgraf
	{\it E-mail address:} {\rm mondalshyam055@gmail.com}
}
\date{\today}
\subjclass{Primary 35S30; Secondary 47G30}
\keywords{Fourier integral operator, Bilinear pseudo-differential operator, Modulation space, Boundedness}
\allowdisplaybreaks
\begin{document}
	
	\begin{abstract}
		In this article, we study   properties of multilinear 	Fourier integral operators on weighted modulation spaces. In particular, using the theory of Gabor frames, we study boundedness of multilinear 	Fourier integral operators on products of weighted modulation spaces. Further, we investigate the periodic multilinear Fourier integral operator. Finally, we study continuity of bilinear pseudo-differential operators on modulation spaces for certain symbol classes, namely $\textbf{SG}$-class.
	\end{abstract}
	\maketitle
	\tableofcontents

	\section{Introduction}\label{sec1}
	This paper deals with the boundedness of a  class of multilinear 	Fourier integral operators on weighted modulation spaces. These operators are defined as follows:  If    $\sigma: \mathbb{R}^d \times  \mathbb{R}^{dr}\to \mathbb{C}$  is a measurable function and $\Phi_{i} : \mathbb{R}^d \times \mathbb{R}^{d} \rightarrow \mathbb{R},$  $1 \leq i \leq r,$ are real-valued phase functions, i.e., is homogeneous of degree 1 in the second variable and satisfies the non-degeneracy condition. By non-degeneracy we mean that the modulus of the determinant of the mixed Hessian does not vanish. The multilinear Fourier integral operator associated with   $\sigma$ is  defined by 
	\begin{equation}\label{definition of multilinear}
		T_{\sigma}(f)(x)=\int_{\mathbb{R}^{dr}} e^{2\pi i [\Phi_{1}(x, \xi_1)+\cdots+\Phi_{r}(x, \xi_r)]}  \sigma(x, \xi) \hat{f_1}(\xi_1) \cdots \hat{f_r}(\xi_r) \;d \xi_1 d \xi_2 \cdots d \xi_r ,
	\end{equation}
	where $x, \xi_1,...,\xi_r \in \mathbb{R}^{d}, f=( f_1,f_2,...,f_r) \in \mathcal{S} (\mathbb{R}^d)^r$, and 
	$$\hat{f_i}(\xi_i) =\int_{\mathbb{R}^d} e^{-2\pi i  y\xi_i}f_i(y) dy,\quad  \xi_{i}\in \mathbb{R}^d, $$
	is the usual Fourier transform of $f_i$. The function $\sigma$  in (\ref{definition of multilinear}) is usually referred as the amplitude of the operator $T_{\sigma}$. On the other hand, if $\sigma: \mathbb{T}^d\times  \mathbb{Z}^{dr}\to \mathbb{C}$  be a measurable function and $\Phi_{i} : \mathbb{T}^d \times \mathbb{Z}^{d} \rightarrow \mathbb R,$  $1 \leq i \leq r,$ are 	   real-valued phase functions    such that   $\Phi_{i}$'s are linear in the second variable for each $i$, then  the periodic multilinear Fourier integral operator defined as 
	$$
		T_{\sigma}(f )(x)=\sum_{k\in  \mathbb{Z}^{dr}} e^{2\pi i [\Phi_{1}(x, k_1)+\cdots+\Phi_{r}(x, k_r)]}  \sigma(x, k) \hat{f_1}(k_1) \cdots \hat{f_r}(k_r), \quad x \in \mathbb{T}^d,
	$$
	where   $ k=(k_1,...,k_r)\in \mathbb{Z}^{dr}, f=( f_1,f_2,...,f_r) \in  C^{\infty} (\mathbb{T}^d)^r$, and $$\hat{f_i}(k_i) = \int_{\mathbb{T}^d} e^{-2 \pi i \eta\cdot k_{i}} f_{i}(\eta), d \eta, \quad k_i\in \mathbb{Z}^d,$$  is the periodic  Fourier transform of $f_i$.
	
	For the Euclidean case, when $r \geq  2$, these operators have been studied by    Rodríguez-López, Rule,  and   Staubach  \cite{RRS}. Particularly, for $r=2$, bilinear Fourier integral  operators have been investigated  by several authors; we refer to \cite{GP, RRS1,HZ} and references therein. 
	If $r = 1$, these quantization formulae can be reduced to the familiar    Fourier integral operator    defined   in the following way:
	\begin{align}\label{FIOs}
		T f(x)=\int a(x, \xi) \widehat{f}(\xi) \mathrm{e}^{i \Phi(x, \xi)} d \xi, \qquad x\in \mathbb{R}^n,
	\end{align}
	where $\hat{f}$ is the Fourier transform of a Schwartz function $f, a$ is the amplitude function,  and $\Phi$ is the phase function. In the literature  of FIOs, authors often considered  the amplitude  function $a$ from the H\"ormander   class  $ S_{\rho, \delta}^m$,  consists of     smooth functions $a \in C^{\infty}\left(\mathbb{R}^d \times \mathbb{R}^d\right)$   satisfying
	$$
	\left|\partial_x^\alpha \partial_{\xi}^\beta a(x, \xi)\right| \leq C_{\alpha, \beta}(1+|\xi|)^{m+\delta|\alpha|-\rho|\beta|},
	$$
	for  all multi-indices $\alpha$ and $\beta$, where  $m\in \mathbb{R}$ and $\rho, \delta\in [0,1].$ The general theory of FIOs  was first developed by Duistermaat and H\"ormander (see \cite{LH,DL}). The theory of   these operators  plays an important role in Fourier analysis and solves various problems arising in partial differential equations \cite{Ruzhansky1}. For instance, in estimating the solutions to a long variety of hyperbolic problems, the mapping properties of FIOs is a   fundamental task, see \cite{Stein}. They were widely employed to study the spectral property of a class of globally elliptic operators, generalizing the harmonic oscillator of the Quantum mechanics
	\cite{Helffer}. The subject of FIOs has been studied by several authors; we refer the reader to \cite{Sogge, Ruzhansky, Fujiwara, Tao} and references therein. Note that when the phase function is given by $\phi(x, \xi) =  x\cdot \xi$, then the   (\ref{FIOs}) reduces to the pseudo-differential operators and Fourier multipliers \cite{Hormander, Ruzhansky1}.  
	
	A fundamental problem in the theory of FIOs is that of classifying the interplay between the properties of a symbol and the properties of its associated Fourier integral operator. Thus one of the natural questions in the theory of FIOs is to find sufficient (nontrivial) conditions on the amplitude function such that the corresponding Fourier integral operator will be bounded on products of certain Banach spaces such as Lebesgue, Sobolev, or Besov spaces.   The local $L^2$ regularity theory for FIOs  with smooth amplitude goes back to the work of Èskin \cite{Eskin}. Further,  a  global $L^2$-boundedness  was investigated by  H\"ormander \cite{Hormander} (see also \cite{Ruzhansky2}). More generally,  Ruzhansky and Sugimoto \cite{Ruzhansky3} studied weighted Sobolev $L^2$-estimates for the FIOs. 
	
	In the direction of $L^p$,  Seeger, Sogge, and Stein   \cite{Sogge}    investigated the local $L^p, 1<p<\infty,$ boundedness of smooth FIOs.   We refer to the book of Sogge  \cite{Sogge1} for an extensive study on the local boundedness of FIOs. A weak- $L^1$ estimate was proved by Tao \cite{Tao}. The global $L^p$ boundedness was established by Cordero, Nicola, and Rodino in \cite{Cordero} when the amplitude function is in the  SG classes. For the general amplitudes from the classes $L^p S_{\rho, \delta}^m$ (i.e,  for rough symbol), FIOs have been considered by Coriasco and Ruzhansky \cite{Coriasco}.

	For bilinear FIOs,    Grafakos and Peloso \cite{GP} proved that if the phases are homogeneous of degree one in the frequency variables with non-degeneracy condition, then under some assumptions on the amplitude functions,  the corresponding FIO  is $L^p\times L^q\to L^r$  bounded for $m<-(n-1)\left(\left|\frac{1}{p}-\frac{1}{2}\right|+\left|\frac{1}{q}-\frac{1}{2}\right|\right)$, $\frac{1}{p}+\frac{1}{q}=\frac{1}{r}$ with $1\leq p,q\leq 2$. Further,   the authors in \cite{Rodriguez1} proved that  result of \cite{GP} could be extended to a global $L^p\times L^q\to L^r$  boundedness  for the full range of exponents $1\leq p,q\leq \infty$ including the endpoint $m=-(n-1)\left(\left|\frac{1}{p}-\frac{1}{2}\right|+\left|\frac{1}{q}-\frac{1}{2}\right|\right)$.  Later, a  bilinear analog of   Seeger–Sogge–Stein theorem, i.e., boundedness theorem to the endpoints of $m$ for bilinear FIOs has been proved by \cite{Rodriguez1}. Recently, the global boundedness of a class of multilinear FIOs was investigated by the authors in  \cite{RRS}, and this one is the only paper to date related to the regularity of multilinear FIOs.

	Motivated by the work on the bilinear pseudo-differential operators on Modulation spaces \cite{Okoudjou}, in this paper, we also investigate multilinear  FIOs on weighted Modulation spaces.    The modulation spaces were first introduced by Feichtinger \cite{fei83, fei97} by imposing integrability conditions on the STFT of tempered distributions. Note that the modulation spaces play a crucial role in the theory of Gabor frames.  Gabor frames provide very efficient representations for a large class of FIOs, mainly to study the boundedness of FIOs on modulation spaces. Moreover, modulation spaces were recently used to formulate and prove boundedness results of linear FIOs \cite{EC&FN&LR}. This is yet another motivation to study the boundedness properties of multilinear FIOs in terms of modulation spaces.

	Particularly,    we prove the boundedness of multilinear  FIOs  on products of modulation spaces.   We show that multilinear  FIOs  corresponding to non-smooth symbols in the Feichtinger algebra are bounded on products of modulation spaces.	      The approach we use here is fundamentally different from the ones previously employed in dealing with the linear case; namely,   we also employ decomposition techniques of functions spaces; however, our first novelty of this paper is that here we use the theory of Gabor frames expansions of tempered distributions in the so-called modulation spaces to prove our boundedness results. We decompose the functions in the modulation spaces into their Gabor expansions and thereby transform the boundedness of the multilinear operator into that of an infinite matrix acting on sequence spaces associated with the modulation spaces. The conditions we impose on the infinite matrix to prove our results turn out to be equivalent to the membership of the corresponding amplitude functions to a particular modulation space. The second novelty is that the results are proved for multilinear and not just bilinear operators. The following subsection gives the main results of this paper.
	\subsection{Main results} 	Let $T_{\sigma}$ be a multilinear Fourier integral operator associated to a symbol $\sigma$. Then    $T_{\sigma}$ coincides with a multilinear integral operator $B_{K}$ with kernel $K.$    The following result is about a connection between the symbol of the multilinear Fourier integral operator and its corresponding integral kernel in terms of modulation spaces.
	\begin{theorem} 
		Let $T_{\sigma}$ be a multilinear Fourier integral operator associated with the  amplitude function  $\sigma$. Then 
		$\sigma(x, y_1, y_2,...,y_r) e^{2\pi i \Phi_{1}(x, y_1)} \cdots e^{2\pi i \Phi_{r}(x, y_r)} \in M_{\Omega_{s}^{B}}^{1}\left(\mathbb{R}^{d(r+1)}\right)$ if and only if $K(x, y_1, y_2,...,y_r)
		\in M_{\Omega_{s}}^{1}\left(\mathbb{R}^{d(r+1)}\right).$
		Moreover, if each $ \Phi_{i}(x, y_i), 1 \leq i \leq r$, are continuous function on $\mathbb{R}^{2d}$, then we have $\sigma \in M_{\Omega_{s}^{B}}^{1}\left(\mathbb{R}^{d(r+1)}\right)$ if and only if $K(x, y_1, y_2,...,y_r)
		\in M_{\Omega_{s}}^{1}\left(\mathbb{R}^{d(r+1)}\right).$
	\end{theorem}
	The following result is about the boundedness of a multilinear integral operator with kernel in the modulation space $M_{\Omega_{s}}^{1}$.
	\begin{theorem} 
		Let $v$ be an $s$-moderate weight, and $1 \leq p_{i}, q_{i}, s_{t}<\infty,$ for $1 \leq i \leq r,$ and $t \in \{1,2\}$, be such that $\frac{1}{p_{1}}+\cdots+\frac{1}{p_{r}}=\frac{1}{s_{1}}$ and $\frac{1}{q_{1}}+\cdots+\frac{1}{q_{r}}=\frac{1}{s_{2}}$. If $K \in M_{\Omega_{s}}^{1}\left(\mathbb{R}^{d(r+1)}\right)$, then the multilinear integral operator $B_{K}$  can be extended as a bounded operator from $M_{v}^{p_{1}, q_{1}}\left(\mathbb{R}^{d}\right) \times \cdots \times M_{v}^{p_{r}, q_{r}}\left(\mathbb{R}^{d}\right)$ into $M_{v}^{s_{1}, s_{2}}\left(\mathbb{R}^{d}\right)$.
	\end{theorem}
	An immediate consequence of the above result provides a sufficient condition on the symbol so that the corresponding FIO  is bounded on products of modulation spaces.
	\begin{theorem}
		Let $v$ be an $s$-moderate weight, and let $1 \leq p_{i}, q_{i}, s_{t}<\infty,$ for $1 \leq i \leq r,$ and $ t \in \{1,2\}$,  be such that $\frac{1}{p_{1}}+\cdots+\frac{1}{p_{r}}=\frac{1}{s_{1}}$ and $\frac{1}{q_{1}}+\cdots+\frac{1}{q_{r}}=\frac{1}{s_{2}}$.  Let $\sigma_{0}(t_1, t_2,...,t_r, t_{r+1})=\sigma(t_1, t_2,...,t_r, t_{r+1}) $ $\times e^{2\pi i \Phi_{1}(t_1,t_2)}\cdots e^{2\pi i \Phi_{r}(t_1, t_{r+1})}$. If $\sigma_{0} \in M_{\Omega_{s}^{B}}^{1}\left(\mathbb{R}^{d(r+1)}\right)$,  then the corresponding  Fourier integral operator $T_{\sigma_{0}}$  can be extended as a bounded operator from $M_{v}^{p_{1}, q_{1}}\left(\mathbb{R}^{d}\right) \times \cdots \times M_{v}^{p_{r}, q_{r}}\left(\mathbb{R}^{d}\right)$ into $M_{v}^{s_{1}, s_{2}}\left(\mathbb{R}^{d}\right)$.
	\end{theorem}
	The following result is about the boundedness of bilinear pseudo-differential operators on modulation spaces for certain symbol classes, namely $\mathbf{SG}$-class, denoted by $\mathbf{SG}^{m_{1}, m_{2}, m_{3}}$. 
	\begin{theorem} 
		For $s_{1}<<0, s_{2}>0$, let $\mu \in \mathcal{M}_{v_{s_{1}, s_{2}}}$.
		Consider a symbol $\sigma$ satisfying
		$$	\left|\partial_{x}^{\alpha} \partial_{\xi}^{\beta} \partial_{\eta}^{\gamma} \sigma(x, \xi, \eta)\right| \leq C_{\alpha, \beta, \gamma}\langle x \rangle^{m_{3}}\langle \xi\rangle^{m_{1}} \langle \eta\rangle^{m_{2}}, \quad|\alpha| \leq 2 N_{3},|\beta| \leq 2 N_{1},|\gamma| \leq 2 N_{2},$$
		
		with $N_{1}>\frac{s_2 +d}{2}$, $N_{2}>\frac{d}{2}$, $\frac{|m_2|+d}{2} < N_{3} < \frac{-s_1-|m_1|-d}{2}$. Then, for every $1 \leq p, q \leq \infty,$ the corresponding pseudo-differential operatr $T_\sigma$ extends to a continuous operator from $\tilde{M}_{\mu}^{p, q} \otimes \tilde{M}_{\mu}^{p, q}$ to $\tilde{M}_{\mu v_{-m_{1}-m_{2}-2 N_{3},-m_3}}^{p, q}$.
	\end{theorem}

	Apart from the introduction, the paper is organized as follows. 
	\begin{itemize}
		\item In Section \ref{sec2},   we recall the definition of modulation spaces and some of their important properties. We also define  Gabor frames and collect some of its properties that we are going to use in this article. 
		
		\item In Section \ref{sec3},   we investigate the relationship between a multilinear integral operator having a kernel
		and a multilinear Fourier integral operator. In particular, we show that the multilinear Fourier integral operator $T_{\sigma}$ coincides with a bilinear integral operator $B_{K}$ having kernel $K$.
		\item In Section \ref{sec4},   we study the boundedness property of the multilinear 	Fourier integral operator on modulation spaces using the theory of Gabor frames.
		
		\item In Section \ref{sec5},    we investigate the periodic multilinear Fourier integral operator on modulation spaces over $\mathbb{T}^d.$
		
		\item In Section \ref{sec6},     we study the boundedness of bilinear pseudo-differential operators on modulation spaces for certain symbol classes, namely $\mathbf{SG}^{m_{1}, m_{2}, m_{3}}$. 
	\end{itemize} 
	

	\section{Preliminaries}\label{sec2}
	In this section we present some basics on the theory of  Modulation space and Gabor Frames.
	
	\subsection{Modulation spaces}
	In this subsection, we briefly recall modulation space over $\mathbb{R}^d$. For a complete background about Modulation spaces, we refer to the reader \cite{fei83, fei97}, and references therein. 
	
	Let $  \mathcal{S}\left(\mathbb{R}^d\right)$ be the Schwartz space of rapidly decreasing functions, with its usual Fr\'{e}chet topology and     the  dual  of $  \mathcal{S}\left(\mathbb{R}^d\right)$, i.e.,  the set of all tempered distributions on $\mathbb{R}^d$ is denoted by $  \mathcal{S}^{\prime}\left(\mathbb{R}^d\right)$.  
	
	For $x, \xi \in \mathbb{R}^d$, let $M_\xi$ and $T_x$  defined as 
	$$T_xf (t) = f (t-x)~~\quad \text{and}~~ \quad M_\xi f (t) = e^{2\pi i \xi \cdot t} f (t),$$
	denote the modulation and translation operators, respectively.  Then, the short-time Fourier transform (STFT)  of a function $f$ with respect to a window function $g \in  \mathcal{S}(\mathbb{R}^d)$ is defined by
	\[ V_g f (x,\xi)=\langle f, M_\xi T_x g \rangle=\int_{\mathbb{R}^d} f(t) \overline{g(t-x)} e^{-2 \pi i \xi \cdot t} dt, \quad (x,\xi) \in \mathbb{R}^d. \]

	For  a strictly positive function $v$ on $\mathbb{R}^{2 d}$, we denote  $L_v^{p, q}=L_v^{p, q}\left(\mathbb{R}^{2 d}\right)$ as the spaces of measurable functions $f$ for which the  following weighted mixed norms
	$$
	\|f\|_{L_v^{p, q}}=\left(\int_{\mathbb{R}^d}\left(\int_{\mathbb{R}^d}|f(x, y)|^p v(x, y)^p d x\right)^{q / p} d y\right)^{1 / q}
	$$
	are finite. When  $p=q$,  the $L_v^{p, p}\left(\mathbb{R}^{2 d}\right)=L_v^p\left(\mathbb{R}^{2 d}\right)$ is the usual  weighted Lebesgue space. Similarly,  $l_{\tilde{v}}^{p, q}\left(\mathbb{Z}^{2 d}\right)$,   denotes the spaces of sequences $a=\left(a_{k l}\right)_{k, l \in \mathbb{Z}^d}$ for which the  following mixed norms
	$$
	\|a\|_{l_{\tilde{v}}^{p, q}}=\left(\sum_{k \in \mathbb{Z}^d}\left(\sum_{l \in \mathbb{Z}^d}\left|a_{k l}\right|^p \tilde{v}(k, l)^p\right)^{q / p}\right)^{1 / q}
	$$
	are finite, where $\tilde{v}(k, l)=v(\alpha k, \beta l)$ for some fixed $\alpha, \beta>0$. If $p=q$, then  $l_{\tilde{v}}^p\left(\mathbb{Z}^{2 d}\right)$  is  the usual weighted sequence spaces.
	\begin{defn}
		A positive, even weight function $\mu \neq 0$ everywhere on $\mathbb{R}^{2 d}$ belongs to $\mathcal{M}_v$ if it satisfies the following condition for some $C>0$:
		$$
		\mu\left(x+y\right) \leq C v\left(x\right) \mu\left(y\right) \quad \forall x,y \in \mathbb{R}^{2 d}.
		$$
	\end{defn}
	Also, associated to every submultiplicative weight, we consider the class of so-called $s$-moderate weights.
	\begin{defn}
		Let $s \geq 0$.  A positive, continuous, and symmetric function $v$ is called an $s$-moderate weight if 
		$$
		v(x+y) \leq C\left(1+|x|^2\right)^{s / 2} v(y), \quad 	\forall x, y \in \mathbb{R}^d, 
		$$
		for some   constant $C>0.$
	\end{defn}
	
	The modulation spaces were introduced by Feichtinger \cite{fei83, fei97}, by imposing integrability conditions on the STFT of tempered distributions.
	Here, we are interested in weighted modulation spaces with respect to the  $s$-moderate weight $v$ defined on $\mathbb{R}^{2d}$. For $s$-moderate weight  function $v$, the mixed Lebesgue space $L^{p,q}_{v}(\mathbb{R}^{2d}), 1\leq p,q \leq \infty,$  is defined by
	$$
	L^{p,q}_{v}(\mathbb{R}^{2d})=\{ f: f \text{ is measurable on } \mathbb{R}^{2d}  ~\text{such that}~ \Vert f \Vert_{L^{p,q}_{v}(\mathbb{R}^{2d})} < \infty\},
	$$
	where $\Vert f \Vert_{L^{p,q}_{v}(\mathbb{R}^{2d})} = \displaystyle \bigg( \int_{\mathbb{R}^d} \bigg( \int_{\mathbb{R}^d} |f(x,\xi)|^p |v(x,\xi)|^p  dx \bigg)^{q/p}  d\xi \bigg)^{1/q},$ with the usual modifications if $p$ and/or $q$ are infinite. When $p=q$, we donote $L^{p,p}_{v}(\mathbb{R}^{2d})=L^{p}_{v}(\mathbb{R}^{2d})$
	\begin{defn}
		Let $v$  be a  $s$-moderate weight function  on $\mathbb{R}^{2d}$,  $g \in \mathcal{S}(\mathbb{R}^d)$ be a  fixed non-zero window function, and $1 \leq p,q \leq \infty$. Then the weighted modulation space $M_v^{p,q}(\mathbb{R}^d )$  consists of all tempered distributions $f \in \mathcal{S'}(\mathbb{R}^d)$ such that $V_g f \in L_v^{p,q}(\mathbb{R}^{2d})$. The norm on $M_v^{p,q}(\mathbb{R}^d )$ is given by
		\begin{eqnarray*}
			\Vert f \Vert_{M_v^{p,q}(\mathbb{R}^d )}
			&=& \Vert V_g f \Vert_{L_v^{p,q}(\mathbb{R}^{2d}  )} \\
			&=& \bigg( \int_{\mathbb{R}^d} \bigg( \int_{\mathbb{R}^d} |V_gf(x,\xi)|^p |v(x,\xi)|^p  dx \bigg)^{q/p}  d\xi \bigg)^{1/q} < \infty,
		\end{eqnarray*}
		with the usual adjustments if $p$ or $q$ is infinite.
	\end{defn}
	If $p=q$, then we write $M_v^p(\mathbb{R}^d )$ instead of $M_v^{p,p}(\mathbb{R}^d )$. When $v=1$ on $\mathbb{R}^d$, then we write $M^{p, q}(\mathbb{R}^d )$ and $M^p(\mathbb{R}^d )$  for   $M_v^{p, q}(\mathbb{R}^d )$ and   $M_v^p(\mathbb{R}^d ),$  respectively.   The definition of $M_v^{p, q}(\mathbb{R}^d )$   is  independent of the choice of $g$ in the sense that each different choice of $g$ defines   equivalent norms on $M_v^{p, q}(\mathbb{R}^d )$. Each weighted modulation space is a Banach space. For $p=q=2$, we have   $M_v^2(\mathbb{R}^d ) =L_v^2(\mathbb{R}^d ).$ 
	If $v(x, \xi)=\left(1+|\xi|^2\right)^{s / 2}$, then $M_v^2(\mathbb{R}^d )=H^s(\mathbb{R}^d)$, the standard Sobolev space. Moreover, if $v(x, \xi)=\left(1+|x|^2+|\xi|^2\right)^{s / 2}$ then $M_v^2(\mathbb{R}^d )=L_s^2(\mathbb{R}^d ) \cap H^s(\mathbb{R}^d )$. Also, note that we can generalize the definition of $s$-moderate weight functions to higher dimensions in the following manner: Let $\omega_s=\left(1+|x|^2+|\xi|^2\right)^{s / 2}$ defined over $\mathbb{R}^{2d}$, we let $\Omega_s=\omega_s \otimes \omega_s \otimes \omega_s$, i.e.,
	$$
	\Omega_s\left(x_1, \xi_1, x_2, \xi_2, x_3, \xi_3\right)=\omega_s\left(x_1, \xi_1\right) \omega_s\left(x_2, \xi_2\right) \omega_s\left(x_3, \xi_3\right),
	$$
	be a weight function defined on $\mathbb{R}^{6d}$. Moreover, if $A$ is an invertible transformation on $\mathbb{R}^{6d}$, we denote by $\Omega_s^A$ the weight function defined on $\mathbb{R}^{6d}$ by 
	\begin{eqnarray}\label{transform defi}
		\Omega_s^A(X)=\Omega_s(A(X)), \quad X \in \mathbb{R}^{6d}.
	\end{eqnarray}
	
	We already saw that for $p=q=2$, we have   $M_v^2(\mathbb{R}^d ) =L_v^2(\mathbb{R}^d ).$ For other $p=q$, the space $M_v^p(\mathbb{R}^d )$ is not $L_v^p(\mathbb{R}^d )$.    In fact for $p=q>2$, the space $M_v^p(\mathbb{R}^d )$ is a superset of $L_v^2(\mathbb{R}^d )$. We have the following inclusion  \[ \mathcal{S}(\mathbb{R}^d) \subset M_v^1(\mathbb{R}^d ) \subset M_v^2(\mathbb{R}^d )=L_v^2(\mathbb{R}^d ) \subset M_v^\infty(\mathbb{R}^d ) \subset \mathcal{S'}(\mathbb{R}^d). \] 	  Particularly, we have $M_v^p(\mathbb{R}^d ) \hookrightarrow L_v^p(\mathbb{R}^d )$ for $1 \leq p \leq 2$, and  $L_v^p(\mathbb{R}^d ) \hookrightarrow M_v^p(\mathbb{R}^d )$ for $2 \leq p \leq \infty$. Let $\tilde{M}_{v}^{p, q} $ denotes the closure of the Schwartz class in $M_v^{p,q}(\mathbb{R}^d ).$ Moreover, the dual of a weighted modulation space is also a weighted modulation space, if $p < \infty$, $q < \infty$, $(M_v^{p, q}(\mathbb{R}^d ))^{'} =M_{\frac{1}{v}}^{p', q'}(\mathbb{R}^d )$, where $p', \; q'$ denote the dual exponents of $p$ and $q$, respectively. 
	
	We note that the modulation space $M_{\omega_s}^1$ is a Banach algebra under both pointwise multiplication and convolution. Moreover, it is    invariant under the Fourier transform.  Modulation space plays also a crutial  role in the theory of Gabor frames as  it serves as a convenient class of windows that generate Gabor frames for the whole class of the modulation spaces.
	In particular, if $s=0$ (equivalently, if $v=\omega_s \equiv 1$ ), then $M^1$ is the Feichtinger algebra. The  functions in $M_{\omega_s}^1$ are not smooth in general.  
	
	Similar to the Euclidean case,    the modulation spaces on $\mathbb{T}^d$ are defined using the STFT on $\mathbb{T}^d \times \mathbb{Z}^d$. Fix a non-zero window $g \in C^\infty\left(\mathbb{T}^d\right)$, and, let $1 \leq p, q \leq \infty$. Then the modulation space $M^{p,q}\left(\mathbb{T}^d\right)$ consists of all tempered distributions $f \in \mathcal{D}^{\prime}\left(\mathbb{T}^{d}\right)$ such that $V_g f \in L^{p,q}\left(\mathbb{T}^d \times \mathbb{Z}^d\right)$. The norm on $M^{p,q}\left(\mathbb{T}^d\right)$ is
	$$
	\|f\|_{M^{p,q}\left(\mathbb{T}^d\right)}=\left\|V_g f\right\|_{L^{p,q}\left(\mathbb{T}^d \times \mathbb{Z}^d\right)}=\left(\int_{\mathbb{T}^d} \left( \sum_{m \in \mathbb{Z}^d}  \left|V_g f(w, m)\right|^p \right)^{^{q / p}}d w\right)^{1 / q}<\infty
	$$
	with the usual adjustments if $p,q$ are infinite. The definition of $M^{p,q}\left(\mathbb{T}^d\right)$ is independent of the choice of $g$ in the sense that each different choice of $g$ defines an equivalent norm on $M^{p,q}\left(\mathbb{T}^d\right)$. Moreover,   all the properties for modulation space over $\mathbb{T}^d$ holds as similar to modulation space over $\mathbb{R}^d$.  For further properties and uses of weighted modulation spaces, we refer to \cite{GK, Okoudjou, fei83}.
	
	\subsection{Gabor frames}
	Fix a function $g \in L^{2}\left(\mathbb{R}^{d}\right)$ and a lattice $\Lambda=\alpha \mathbb{Z}^{d} \times \beta \mathbb{Z}^{d}$, for $\alpha, \beta>0$. For $(m, n) \in \Lambda$, define $g_{m, n}:=M_{n} T_{m} g$. The set of time-frequency shifts $\mathcal{G}(g\otimes g, \alpha, \beta)=\left\{(g_{m, n},g_{m_{0},n_{0}});(m, n),(m_0,n_0) \in \Lambda\right\}$ is called bilinear Gabor system. We define the coefficient operator $C_{g}$, associated to $\mathcal{G}(g\otimes g, \alpha, \beta)$, which maps functions to sequences as follows:
	$$
	\begin{aligned}
		\left(C_{g} (f_1,f_2)\right)_{m, n,m_0,n_0}&=\left((C_{g}^{\alpha, \beta} f_1)_{m, n} , (C_{g}^{\alpha, \beta} f_2)_{m_0, n_0}\right)\\
		&:=\left(\left\langle f_1, g_{m, n}\right\rangle, \left\langle f_2, g_{m_0, n_0}\right\rangle\right), \quad(m, n),(m_0, n_0) \in \Lambda.
	\end{aligned}
	$$
	Also, define the synthesis operator as follow:
	$$
	D_{g} (c,e)=D_{g}^{\alpha, \beta}(c,e)=\sum_{(m, n),(m_0, n_0) \in \Lambda} c_{m, n} e_{m_0, n_0} (M_{n} T_{m} g, M_{n_0} T_{m_0} g),
	$$
	where $c=\left\{c_{m, n}\right\}_{(m, n) \in \Lambda}, e=\left\{e_{m_0, n_0}\right\}_{(m_0, n_0) \in \Lambda}.$ The Gabor frame operator corresponding to the coefficient operator and synthesis operator is defined by
	$$
	S_{g}  (f_1,f_2)=S_{g}^{\alpha, \beta}  (f_1,f_2):=D_{g} C_{g}  (f_1,f_2)=\sum_{(m, n),(m_0, n_0) \in \Lambda}\left\langle f_1, g_{m, n}\right\rangle \overline{\left\langle f_2, g_{m_0, n_0}\right\rangle} ( g_{m, n}, g_{m_0, n_0}) .
	$$
	The set $\mathcal{G}(g \otimes g, \alpha, \beta)$ is called a Gabor frame for the Hilbert space $L^{2}\left(\mathbb{R}^{d}\right) \otimes L^{2}\left(\mathbb{R}^{d}\right)$ if $S_{g}$ is a bounded and invertible operator on $L^{2}\left(\mathbb{R}^{d}\right) \otimes L^{2}\left(\mathbb{R}^{d}\right)$. If $\mathcal{G}(g\otimes g, \alpha, \beta)$ is a Gabor frame for $L^{2}\left(\mathbb{R}^{d}\right) \otimes L^{2}\left(\mathbb{R}^{d}\right)$, then the so-called dual window $\gamma=S_{g}^{-1} g$ is well defined and the set $\mathcal{G}(\gamma \otimes \gamma, \alpha, \beta)$ is called the canonical dual frame of $\mathcal{G}(g \otimes g, \alpha, \beta))$. Every $(f_1,f_2) \in L^{2}\left(\mathbb{R}^{d}\right) \otimes L^{2}\left(\mathbb{R}^{d}\right)$ possesses the non-orthogonal frame expansion
	$$
(f_1,f_2)=\sum_{(m, n),(m_0, n_0) \in \Lambda}\left\langle f_1, g_{m, n}\right\rangle \overline{\left\langle f_2, g_{m_0, n_0}\right\rangle} ( \gamma_{m, n}, \gamma_{m_0, n_0}),\\
	$$
	and
	$$
	(f_1,f_2)=\sum_{(m, n),(m_0, n_0) \in \Lambda}\left\langle f_1, \gamma_{m, n}\right\rangle \overline{\left\langle f_2, \gamma_{m_0, n_0}\right\rangle} (g_{m, n},g_{m_0, n_0})
	$$
	with both sums converge unconditionally in $L^{2}\left(\mathbb{R}^{d}\right) \otimes L^{2}\left(\mathbb{R}^{d}\right)$, and norm equivalence:
	$$
	\|(f_1,f_2)\|_{L^{2} \otimes L^{2}} \asymp\left\|C_{g} (f_1,f_2)\right\|_{\ell^{2} \otimes \ell^{2}} \asymp\left\|C_{\gamma} (f_1,f_2)\right\|_{\ell^{2} \otimes \ell^{2}} .
	$$
	In particular, if $\gamma=g$ and $\|g\|_{L^{2}}=1$, then the frame is called normalized tight bilinear Gabor frame. The following theorem gives us the boundedness of analysis operator and synthesis operator on modulation spaces and sequence spaces, respectively.
	\begin{theorem}\label{gabor frame bdd}
		Let $\mu \in \mathcal{M}_{v},$ and $\mathcal{G}(g \otimes g, \alpha, \beta)$ be a normalized tight bilinear Gabor frame for $L^{2}\left(\mathbb{R}^{d}\right) \otimes L^{2}\left(\mathbb{R}^{d}\right)$, with lattice $\Lambda=\alpha \mathbb{Z}^{d} \times \beta \mathbb{Z}^{d}$, where $g \in M_{v}^{1}$. Define $\tilde{\mu}=\mu_{\left.\right|_{\Lambda}}$ and let $\tilde{\ell}_{\tilde{\mu}}^{p, q}$ denotes the closure of the space of eventually zero sequences in $\ell_{\tilde{\mu}}^{p, q}$. Then for every $1 \leq p, q \leq \infty$ the operator $C_{g}$ is continuous from $\tilde{M}_{\mu}^{p, q} \otimes \tilde{M}_{\mu}^{p, q}$ into $\tilde{\ell}_{\tilde{\mu}}^{p, q} \otimes \tilde{\ell}_{\tilde{\mu}}^{p, q}$, whereas the operator $D_{g}$ is continuous from $\tilde{\ell}_{\tilde{\mu}}^{p, q} \otimes \tilde{\ell}_{\tilde{\mu}}^{p, q}$ into $\tilde{M}_{\mu}^{p, q} \otimes \tilde{M}_{\mu}^{p, q}.$
	\end{theorem}
	Let us now define the following symbol class which is just the multilinear version of H\"ormander symbol class on $\mathbb{R}^d.$
	\begin{defn}\label{symbol class on R^n}
		Let $0 \leq \rho, \delta \leq 1$. We say a symbol $\sigma(x, \xi_1, \xi_2,...,\xi_r) \in S_{\rho,\delta}^{m_1,m_2,...,m_r}\underbrace{\left(\mathbb{R}^d \times \cdots \times \mathbb{R}^d\right)}_{r+1-times},$ where $m_1,m_2,...,m_r \in \mathbb{R}$, if for all multi-indices $\alpha, \beta_1, \beta_2,...,\beta_r$,  there exists a positive constant $C_{\alpha, \beta_1, \beta_2,...,\beta_r}$ such that
		$$|(\partial_{x}^\alpha \partial_{\xi_{1}}^{\beta_{1}} \cdots \partial_{\xi_{r}}^{\beta_{r}}\sigma)(x, \xi_1, \xi_2,...,\xi_r) | \leq C_{\beta_1, \beta_2,...,\beta_r} \prod_{k+1}^{r}\left(\langle \xi_{k} \rangle^{m_{k}-\rho|\beta_k|+\alpha |\delta|}\right) ,\quad  x,\xi_1, \xi_2,...,\xi_r \in \mathbb{R}^d,$$
		and the corresponding Fourier integral operator is defined by \eqref{definition of multilinear}.
	\end{defn}
	\section{Multilinear Fourier integral operators on $\mathbb{R}^d$}\label{sec3}
	In this section,  we study the multilinear Fourier integral operator and find a relation between multilinear integral operators.  We start with the following well defined result. 
	\begin{prop}\label{operator well defined}
		Let $\sigma(x, \xi_1, \xi_2,...,\xi_r) \in S^{m_1,m_2,...,m_r} (\underbrace{\mathcal{S}(\mathbb{R}^{d}) \times \cdots \times  \mathcal{S}(\mathbb{R}^{d})}_{(r+1)-times} ) ,$ where $m_1,m_2,...,m_r \in \mathbb{R}$, and assume that for every multi-index $\alpha$, there exist positive constants $C_1,C_2,...,C_r$ satisfying the following conditions: $$\left|\partial_{x}^{\alpha} e^{i \Phi_{l}(x,\xi_l)}\right| \leq C_l,\qquad 1 \leq l \leq r,$$
		where $\Phi_{l}$'s are real-valued phase functions. 		 
		Then the corresponding multilinear Fourier integral operator $T_{\sigma}$ (defined in  \eqref{definition of multilinear})     is well defined.
	\end{prop}
	\begin{proof}
		The   multilinear Fourier integral operator associated with the amplitude   $\sigma$ is  defined by 
		\begin{equation*}
			T_{\sigma}(f)(x)=\int_{\mathbb{R}^{dr}} e^{2\pi i [\Phi_{1}(x, \xi_1)+\cdots+\Phi_{r}(x, \xi_r)]}  \sigma(x, \xi) \hat{f_1}(\xi_1) \cdots \hat{f_r}(\xi_r) \;d \xi_1 d \xi_2 \cdots d \xi_r ,
		\end{equation*}
		where $x\in \mathbb{R}^d, \xi=(\xi_1,...,\xi_r)\in \mathbb{R}^{dr}, f=( f_1,f_2,...,f_r) \in \mathcal{S} (\mathbb{R}^d)^r$. This also  can be wrtiten as 
		$$
		T_{\sigma}(f)(x)=\int_{\mathbb{R}^{d}} e^{2\pi i \Phi_{1}(x, \xi_1)} \sigma_{1}(x, \xi_1) \hat{f_1}(\xi_1) d \xi_1 ,
		$$
		where 
		\begin{align}\label{symbol}
			\sigma_1(x,\xi_1)=\int_{\mathbb{R}^{d(r-1)}} e^{2\pi i [\Phi_{2}(x, \xi_2)+\cdots+\Phi_{r}(x, \xi_r)]}  \sigma(x, \xi ) \hat{f_2}(\xi_2) \cdots \hat{f_r}(\xi_r)\;d \xi_2 \cdots d \xi_r,\end{align}
		with $\xi=(\xi_1,...,\xi_r)\in \mathbb{R}^{dr}.$ Using Lemma 1.1 of \cite{shubin-book},   substituting  $\left(L_{1}\right)^{k_{1}}e^{2\pi i \Phi_{1}}$ in place of   $e^{2\pi i \Phi_{1}}$ and then integrating by parts $k_{1}$ times, we get 
		$$
		T_{\sigma}(f )(x)=\int_{\mathbb{R}^{d}} e^{2\pi i \Phi_{1}(x, \xi_1)} \left(L_{1}\right)^{k_{1}}\left(\sigma_{1}(x, \xi_1) \hat{f_1}(\xi_1)\right) d \xi_1 ,
		$$
		where
		$$
		L_1 u(x, y)=-\sum_{j=1}^N \frac{\partial}{\partial y_j}\left(a_j u\right)-\sum_{k=1}^n \frac{\partial}{\partial x_k}\left(b_k u\right)+c u,
		$$
		such that $a_j(x, y) \in S^0\left(\mathbb{R}^d \times \mathbb{R}^d\right), b_k(x, y) \in S^{-1}\left(\mathbb{R}^d \times \mathbb{R}^d\right), $ and $c(x,y) \in S^{-1}\left(\mathbb{R}^d \times \mathbb{R}^d\right)$.
		For a moment, we assume that $\sigma_1(x,\xi_1) \in S^{m_1}\left(\mathbb{R}^d \times \mathbb{R}^d\right).$ Then, for every $k_1 \in \mathbb{N},$ we have $$\left(L_{1}\right)^{k_{1}}(\sigma_{1}(x, \xi_1) \hat{f_1}(\xi_1)) \in S^{m_1-k_1\rho}\left(\mathbb{R}^d \times \mathbb{R}^d\right).$$ Now, choose $k_1$ such that $m_1-k_1\rho < -d.$ This gurantees us that the required integral is absolutely convergent. Now it remains to    prove that $\sigma_1(x,\xi_1) \in S^{m_1}\left(\mathbb{R}^d \times \mathbb{R}^d\right).$  From (\ref{symbol}), we have 
		$$
		\begin{aligned}
			\sigma_1(x,\xi_1)&=\int_{\mathbb{R}^{d(r-1)}} e^{2\pi i [\Phi_{2}(x, \xi_2)+\cdots+\Phi_{r}(x, \xi_r)]}  \sigma(x, \xi_1, \xi_2,...,\xi_r) \hat{f_2}(\xi_2) \cdots \hat{f_r}(\xi_r)d \xi_2 \cdots d \xi_r\\
			&=\int_{\mathbb{R}^{d}} e^{2\pi i \Phi_{2}(x, \xi_2)} \sigma_2(x,\xi_1,\xi_2) \hat{f_2}(\xi_2) d \xi_2 ,\\
		\end{aligned}
		$$
		where 
		$$
		\sigma_2(x,\xi_1,\xi_2)=\int_{\mathbb{R}^{d(r-2)}} e^{2\pi i [\Phi_{3}(x, \xi_3)+\cdots+\Phi_{r}(x, \xi_r)]}  \sigma(x, \xi_1, \xi_2,...,\xi_r) \hat{f_3}(\xi_3) \cdots \hat{f_r}(\xi_r)d \xi_3 \cdots d \xi_r.
		$$
		Again, by Lemma 1.1 of \cite{shubin-book},  one can substitute $\left(L_{1}\right)^{k_{2}}e^{2\pi i \Phi_{2}}$ instead of $e^{2\pi i \Phi_{2}}$ and integrate by parts $k_{2}$ times gives 
		$$
		\sigma_1(x,\xi_1)=\int_{\mathbb{R}^{d}} e^{2\pi i \Phi_{2}(x, \xi_2)} \left(L_{1}\right)^{k_{2}}\left(\sigma_{2}(x, \xi_1,\xi_2) \hat{f_2}(\xi_2)\right) d \xi_2 .
		$$
		For a moment, we assume that $\sigma_2(x,\xi_1,\xi_2) \in S^{m_1,m_2}\left(\mathbb{R}^d \times \mathbb{R}^d \times \mathbb{R}^d\right).$ Then, for every $k_2 \in \mathbb{N},$ we have $\left(L_{1}\right)^{k_{2}}(\sigma_{2}(x, \xi_1,\xi_2) \hat{f_2}(\xi_2)) \in S^{m_1,m_2-k_2\rho}\left(\mathbb{R}^d \times \mathbb{R}^d \times \mathbb{R}^d\right).$ Now, choose $k_2$ such that $m_2-k_2\rho < -d.$ Then, using Leibniz’s formula, given condition and the fact  that $\left(L_{1}\right)^{k_{2}}(\sigma(x, \xi_1,\xi_2) \hat{f_2}(\xi_2)) \in S^{m_1,m_2-k_2\rho},$ we obatin $\sigma_1(x,\xi_1) \in S^{m_1}\left(\mathbb{R}^d \times \mathbb{R}^d\right).$ Now, to prove $\sigma_2(x,\xi_1,\xi_2) \in S^{m_1,m_2}\left(\mathbb{R}^d \times \mathbb{R}^d \times \mathbb{R}^d\right),$ we again use the similar techniques developed as above. At the end, we obtained $\sigma_{r}(x, \xi_1, \xi_2,...,\xi_{r-1}) \in S^{m_1,m_2,...,m_{r-1}}\underbrace{\left(\mathbb{R}^d \times \cdots \times \mathbb{R}^d\right)}_{r-times}$, where
		$$
		\sigma_{r}(x, \xi_1, \xi_2,...,\xi_{r-1}) = \int_{\mathbb{R}^{d}} e^{2\pi i \Phi_{r}(x, \xi_r)} \sigma(x, \xi_1,\xi_2,...,\xi_r) \hat{f_r}(\xi_r) d \xi_r,
		$$
		and this completes the proof.
	\end{proof}
	
	We define multilinear integral operator associated with kernel $K \in \mathcal{S}^{\prime}\left(\mathbb{R}^{d(r+1)}\right)$, which is a mapping $B_{K}$ from $\underbrace{\mathcal{S}\left(\mathbb{R}^{d}\right) \times \cdots \times  \mathcal{S}\left(\mathbb{R}^{d}\right)}_{r-times}$ to $\mathcal{S}^{\prime}\left(\mathbb{R}^{d}\right)$ by
	\begin{equation}\label{definition of B_{K}}
		B_{K}(f_1,f_2,...,f_r )(x)=\int_{\mathbb{R}^{dr}} K(x, y_1,...,y_r) f_1(y_1) \cdots f_r(y_r) d y_1 \cdots d y_r,
	\end{equation}
	where $f_1,f_2,...,f_r \in \mathcal{S} (\mathbb{R}^d)$. The next result  establishes the relationship between a multilinear integral operator
	and a multilinear Fourier integral operator.
	\begin{prop}\label{kernel expension}
		Let $T_{\sigma}$ be a multilinear Fourier integral operator associated to a symbol $\sigma$. Then $T_{\sigma}$ coincides with a multilinear integral operator, $B_{K}$ with kernel
		$$K(x, y)=\mathcal{F}_{r+1} \mathcal{F}_{r} \cdots \mathcal{F}_{3} \mathcal{F}_{2}\left(\sigma(x, y ) e^{2\pi i \Phi_{1}(x, y_1)} \cdots e^{2\pi i \Phi_{r}(x, y_r)}\right),$$
		where $y=(y_1, y_2,...,y_r)$ and $ \mathcal{F}_{j}$ denotes the Fourier transform in the $j^{th}$ variable.
	\end{prop}
	\begin{proof}
		For $f=( f_1,f_2,...,f_r)$ with $f_i \in  \mathcal{S}(\mathbb{R}^d), 1\leq i\leq r$, we have
		\begin{align*}
			T_{\sigma}(f )(x)
			&=\int_{\mathbb{R}^{dr}} e^{2\pi i [\Phi_{1}(x, \xi_1)+ \Phi_{2}(x, \xi_2)+\cdots+\Phi_{r}(x, \xi_r)]}  \sigma(x, \xi ) \hat{f_1}(\xi_1) \hat{f_2}(\xi_2) \cdots \hat{f_r}(\xi_r) d \xi_1 d \xi_2 \cdots d \xi_r \\
			&=\int_{\mathbb{R}^{2dr}} e^{2\pi i [\Phi_{1}(x, \xi_1)- \xi_1 \cdot y_1+\cdots+\Phi_{r}(x, \xi_r)- \xi_r\cdot y_r]}  \sigma(x, \xi )   f(y_1) \cdots f(y_r) d \xi_1\cdots d \xi_r d y_1 \cdots d y_r\\
			&=\int_{\mathbb{R}^{dr}} K(x, y_1,\cdots,y_r) f_1(y_1) \cdots f_r(y_r) d y_1 \cdots d y_r\\
			&=B_{K}(f  )(x),
		\end{align*}
		where
		\begin{align*}
			K(x, y_1,\cdots,y_r) &= \int_{\mathbb{R}^{dr}} e^{2\pi i [\Phi_{1}(x, \xi_1)- \xi_1 \cdot y_1+\cdots+\Phi_{r}(x, \xi_r)- \xi_r\cdot y_r]}   \sigma(x, \xi)  d \xi_1\cdots d \xi_r \\
			&=\mathcal{F}_{r+1} \mathcal{F}_{r} \cdots \mathcal{F}_{3} \mathcal{F}_{2}\left(\sigma(x, y_1, y_2,...,y_r ) e^{2\pi i \Phi_{1}(x, y_1)} \cdots e^{2\pi i \Phi_{r}(x, y_r)}\right).\qedhere
		\end{align*}
	\end{proof} 
	Let $u=\left(u_1, u_2,..., u_{r+1}\right), v=\left(v_1, v_2,..., v_{r+1}\right) \in \mathbb{R}^{d(r+1)}$. We define an invertible linear transformation on $\mathbb{R}^{2d(r+1)}$ by $$A(u, v)=\left((u_1,-v_2,...,-v_{r+1}),(v_1,u_2,...,u_{r+1})\right).$$ Let $B$ denotes the inverse of $A$. In the next theorem, we establish a connection between the symbol of the multilinear Fourier intergal operator and its corresponding integral kernel defined in Proposition \ref{kernel expension}, in terms of modulation spaces.
	\begin{theorem}\label{relation of kernel and symbol}	Let $T_{\sigma}$ be a multilinear Fourier integral operator associated  with  the  amplitude  function  $\sigma$. Then 
		$\sigma(x, y_1, y_2,...,y_r) e^{2\pi i \Phi_{1}(x, y_1)} \cdots e^{2\pi i \Phi_{r}(x, y_r)} \in M_{\Omega_{s}^{B}}^{1}\left(\mathbb{R}^{d(r+1)}\right)$ if and only if $K(x, y_1, y_2,...,y_r)
		\in M_{\Omega_{s}}^{1}\left(\mathbb{R}^{d(r+1)}\right),$ where $K$ is same as in Proposition \ref{kernel expension}.
		Moreover, if each $ \Phi_{i}(x, y_i), 1 \leq i \leq r$, are continuous function on $\mathbb{R}^{2d}$, then we have $\sigma \in M_{\Omega_{s}^{B}}^{1}\left(\mathbb{R}^{d(r+1)}\right)$ if and only if $K(x, y_1, y_2,...,y_r)
		\in M_{\Omega_{s}}^{1}\left(\mathbb{R}^{d(r+1)}\right).$
	\end{theorem}
	\begin{proof}
		Let $G \in S\left(\mathbb{R}^{d(r+1)}\right).$ For $u=\left(u_{1},u_{2},...,u_{r},u_{r+1}\right), v=\left(v_{1},v_{2},..., v_{r}, v_{r+1}\right), \text { and } t=\left(t_{1},t_{2},..., t_{r}, t_{r+1}\right)$ $ \in \mathbb{R}^{d(r+1)},$ we obtain
		\begin{align*}
			V_{G} K(u, v)
			&=\int_{\mathbb{R}^{d(r+1)}} e^{-2 \pi i t \cdot v} \overline{G(t-u)} K(t) \;d t \\
			& =\int_{\mathbb{R}^{d(r+1)}} \mathcal{F}_{r+1} \mathcal{F}_{r} \cdots \mathcal{F}_{3} \mathcal{F}_{2}\left(\sigma(t_1, t_2,...,t_r, t_{r+1}) e^{2\pi i \Phi_{1}(t_1,t_2)} \cdots e^{2\pi i \Phi_{r}(t_1, t_{r+1})}\right) \\
			& \qquad \times e^{-2 \pi i\left(t_{1}\cdot v_{1}+\cdots+t_{r+1}\cdot v_{r+1}\right)} \overline{G\left(t_{1}-u_{1},..., t_{r+1}-u_{r+1}\right)} d t_{1} \cdots d t_{r+1}\\
			& =\int_{\mathbb{R}^{d(r+1)}} \mathcal{F}_{r+1} \mathcal{F}_{r} \cdots \mathcal{F}_{3} \mathcal{F}_{2}\left(\sigma_{0}(t_1, t_2,...,t_r, t_{r+1}) \right) e^{-2 \pi i\left(t_{1}\cdot v_{1}+\cdots+t_{r+1}\cdot v_{r+1}\right)}\\
			& \qquad \times \overline{G\left(t_{1}-u_{1},..., t_{r+1}-u_{r+1}\right)} d t_{1} \cdots d t_{r+1},
		\end{align*}
		where $\sigma_{0}(t_1, t_2,...,t_r, t_{r+1})=\sigma(t_1, t_2,...,t_r, t_{r+1}) e^{2\pi i \Phi_{1}(t_1,t_2)} \cdots e^{2\pi i \Phi_{r}(t_1, t_{r+1})}.$ Therefore 
		\begin{align*}
			V_{G} K(u, v)
			& =\int_{\mathbb{R}^{d(r+1)}} \mathcal{F}_{1}^{-1} \widehat{\sigma_{0}}(t_1, t_2,...,t_r, t_{r+1})e^{-2 \pi i\left(t_{1}\cdot v_{1}+\cdots+t_{r+1}\cdot v_{r+1}\right)} \\&\qquad \times\overline{G\left(t_{1}-u_{1},..., t_{r+1}-u_{r+1}\right)} \;d t_{1} \cdots d t_{r+1}\\
			& =\int_{\mathbb{R}^{d(r+2)}}\widehat{\sigma_{0}}(p, t_2,...,t_r, t_{r+1}) e^{2 \pi i p \cdot t_1}e^{-2 \pi i\left(t_{1}\cdot v_{1}+\cdots+t_{r+1}\cdot v_{r+1}\right)} \\&\qquad \times \overline{G\left(t_{1}-u_{1},..., t_{r+1}-u_{r+1}\right)} \;d t_{1} \cdots d t_{r+1} \;d p\\
			& =\int_{\mathbb{R}^{d(r+2)}}\widehat{\sigma_{0}}(p, t_2,...,t_r, t_{r+1}) e^{2 \pi i (u_1-t_1) \cdot (p-v_1)}e^{-2 \pi i\left(t_{2}\cdot v_{2}+\cdots+t_{r+1}\cdot v_{r+1}\right)}\\
			&\qquad \times \overline{G\left(-t_{1},t_2-u_2,..., t_{r+1}-u_{r+1}\right)} \;d t_{1} \cdots d t_{r+1} \;d p\\
			& =\int_{\mathbb{R}^{d(r+1)}}\widehat{\sigma_{0}}(p, t_2,...,t_r, t_{r+1}) e^{2 \pi i u_1 \cdot (p-v_1)}e^{-2 \pi i\left(t_{2}\cdot v_{2}+\cdots+t_{r+1}\cdot v_{r+1}\right)}\\
			&\qquad \times \overline{\mathcal{F}_{1} G\left(p-v_1,t_{2}-u_{2},..., t_{r+1}-u_{r+1}\right)} \;d t_{2} \cdots d t_{r+1} \;d p\\
			& =e^{-2 \pi i u_1 \cdot v_1}\int_{\mathbb{R}^{d(r+1)}}\widehat{\sigma_{0}}(p, t_2,...,t_r, t_{r+1}) e^{-2 \pi i (p,t_2,...,t_{r+1})\cdot (-u_1,v_2,...,v_{r+1})}\\
			&\qquad \times \overline{\mathcal{F}_{1} G\left((p,t_2,...,t_{r+1})-(v_1,u_2,...,u_{r+1})\right)} \;d t_{2} \cdots d t_{r+1} \;d p.
		\end{align*}
		Let $H=\mathcal{F}_{1} G$. We know that $\left|V_{g} f(x, y)\right|=\left|V_{\check{g}} \check{f}(-y, x)\right|,$ whenever the STFT can be defined, so we have
		\begin{align*}
			\left|V_{G} K(u, v)\right|&=\left|V_{H} \widehat{\sigma_{0}}((v_1,u_2,...,u_{r+1}),(-u_1,v_2,...,v_{r+1}))\right| \\
			&=\left|V_{\check{H}}\sigma_{0}((u_1,-v_2,...,-v_{r+1}),(v_1,u_2,...,u_{r+1}))\right|=\left|V_{\check{H}} \sigma_{0}(A(u,v))\right|.
		\end{align*}
		Therefore, by relation \eqref{transform defi}, we get
		\begin{align*}
			\int_{\mathbb{R}^{d(r+1)}} \int_{\mathbb{R}^{d(r+1)}}\left|V_G K(u, v)\right| \Omega_s(u, v)\; d u d v& =\int_{\mathbb{R}^{d(r+1)}} \int_{\mathbb{R}^{d(r+1)}}\left|V_{\check{H}} \sigma_{0}(A(u, v))\right| \Omega_s(u, v) \;d u d v \\
			& =\int_{\mathbb{R}^{d(r+1)}} \int_{\mathbb{R}^{d(r+1)}} \mid V_{\check{H}} \sigma_{0}(u, v) \mid \Omega_s^B(u, v) \;d u d v .
		\end{align*}
		This concludes that $\sigma(x, y_1, y_2,...,y_r) e^{2\pi i \Phi_{1}(x, y_1)} \cdots e^{2\pi i \Phi_{r}(x, y_r)} \in M_{\Omega_{s}^{B}}^{1}\left(\mathbb{R}^{d(r+1)}\right)$ if and only if $$K(x, y)=\mathcal{F}_{r+1} \mathcal{F}_{r} \cdots \mathcal{F}_{2}\left(\sigma(x, y) e^{2\pi i \Phi_{1}(x, y_1)} \cdots e^{2\pi i \Phi_{r}(x, y_r)}\right)\in M_{\Omega_{s}}^{1}\left(\mathbb{R}^{d(r+1)}\right),$$
		where $y=(y_1, y_2,...,y_r).$
		Also, we know that if $f\in M_{w}^{1}$ for any weight $w$ and $\Phi(x, y)$ is a polynoimal type function defined on $\mathbb{R}^{2n}$, then $e^{2\pi i \Phi}\cdot f\in M_{w}^{1}.$ Now, by using the facts that continuous function on $\mathbb{R}^{2n}$ can be approximated by polynomial type functions and the modulation space $M_{w}^{1}$ is a banach space, we get $\sigma \in M_{\Omega_{s}^{B}}^{1}\left(\mathbb{R}^{d(r+1)}\right)$ if and only if $K(x, y_1, y_2,...,y_r)
		\in M_{\Omega_{s}}^{1}\left(\mathbb{R}^{d(r+1)}\right),$ and the proof is complete.
	\end{proof}
	\section{Boundedness of multilinear Fourier integral operators}\label{sec4}
	This section is devoted to obtain the boundedness of multilinear Fourier integral operators on the weighted modulation spaces.

	Consider $\phi \in \mathcal{S}\left(\mathbb{R}^{d}\right)$ that generates a Gabor frame for $L^{2}$ with (canonical) dual $\gamma \in \mathcal{S}\left(\mathbb{R}^{d}\right)$. Then for  $f_1,f_2,...,f_r,h \in \mathcal{S} (\mathbb{R}^d)$,  using the  Gabor series expansion  for multilinear case, we have 
	
	\begin{align}
		\begin{cases}
			&f_{k}=\sum_{m_{k}, n_{k}\in \mathbb{Z}^{d}}\left\langle f_{k}, M_{\beta n_{k}} T_{\alpha m_{k}} \gamma\right\rangle M_{\beta n_{k}} T_{\alpha m_{k}} \phi, \quad 1 \leq k \leq r,\vspace{.4cm}\\ 
			&h=\sum_{i, j\in \mathbb{Z}^{d}}\left\langle h, M_{\beta j} T_{\alpha i} \gamma\right\rangle M_{\beta j} T_{\alpha i} \phi.
		\end{cases}	 
	\end{align}
	
	\noindent From the identity  \eqref{definition of B_{K}} and the above Gabor series expension, we obtain
	$$
	\begin{aligned}
		&	\langle B_{K}(f_1,f_2,...,f_r ),h\rangle\\
		&=\int_{\mathbb{R}^{d(r+1)}} K(x, y_1,...,y_r) \prod_{k=1}^{r} \bigg(\sum_{m_{k}, n_{k} \in \mathbb{Z}^{d}}\left\langle f_{k}, M_{\beta n_{k}} T_{\alpha m_{k}} \gamma\right\rangle M_{\beta n} T_{\alpha m} \phi(y_k)\bigg)\\
		&\qquad  \times \overline{\sum_{i, j \in \mathbb{Z}^{d}}\left\langle h, M_{\beta j} T_{\alpha i} \gamma\right\rangle M_{\beta j} T_{\alpha i} \phi(x)} \;d x d y_1 \cdots d y_r\\
		&= \sum_{i, j} \sum_{m_{1}, n_{1}} \cdots \sum_{m_{r}, n_{r}}\overline{\left\langle h, M_{\beta j} T_{\alpha i} \gamma\right\rangle} \left\langle f_{1}, M_{\beta n_{1}} T_{\alpha m_{1}} \gamma\right\rangle \cdots  \left\langle f_{r}, M_{\beta n_{r}} T_{\alpha m_{r}} \gamma\right\rangle \\
		& \qquad  \times \int_{\mathbb{R}^{d(r+1)}} K(x, y_1,...,y_r) \overline{M_{\beta j} T_{\alpha i} \phi(x)} M_{\beta n_1} T_{\alpha m_1} \phi(y_1) \cdots M_{\beta n_r} T_{\alpha m_r} \phi(y_r) \;
		d x d y_1 \cdots d y_r \\
		&= \sum_{i, j} \sum_{m_{1}, n_{1}} \cdots \sum_{m_{r}, n_{r}}\overline{\left\langle h, M_{\beta j} T_{\alpha i} \gamma\right\rangle} \left\langle f_{1}, M_{\beta n_{1}} T_{\alpha m_{1}} \gamma\right\rangle \cdots  \left\langle f_{r}, M_{\beta n_{r}} T_{\alpha m_{r}} \gamma\right\rangle \\
		&\qquad \times\left\langle B_{K}\left(M_{\beta n_{1}} T_{\alpha m_{1}} \phi,...,M_{\beta n_{r}} T_{\alpha m_{r}} \phi\right), M_{\beta j} T_{\alpha i} \phi\right\rangle .
	\end{aligned}
	$$
	Here, the exchange of the integrals and summations above is justified, since, $f_1,...,f_r, h \in \mathcal{S}$ have absolutely summable Gabor coefficients. Moreover, $K \in \mathcal{S}^{\prime}\left(\mathbb{R}^{d(r+1)}\right)=\bigcup_{s \geq 0} M_{1 / \omega_{s}}^{\infty}$ (cf. \cite[Proposition 11.3.1]{folland}) and $\phi \in \mathcal{S}$ imply that the integral in the second equality is uniformly bounded with respect to $i, j, m_1,n_1,...,m_r,n_r \in \mathbb{Z}^{d}$. Therefore, to study the boundedness of $B_{K}$ on the products of modulation spaces, it suffices to analyze the boundedness of the matrix $B=\left(b_{i j,m_1 n_1,...,m_r n_r}\right)$ defined by
	\begin{equation}\label{infinite matrix form}
		b_{i j,m_1 n_1,...,m_r n_r}=\left\langle B_{K}\left(M_{\beta n_{1}} T_{\alpha m_{1}} \phi,...,M_{\beta n_{r}} T_{\alpha m_{r}} \phi\right), M_{\beta j} T_{\alpha i} \phi\right\rangle .
	\end{equation}
	on products of appropriate sequence spaces.
	
	For an infinite matrix $\left(a_{i j,m_1 n_1,...,m_r n_r}\right)$, let $\mathcal{O}$ denote the multilinear operator associated to it.  {The next theorem is the generalization of  \cite[Theorem 2]{Okoudjou} and can be proved in similar lines.}
	\begin{theorem}\label{bdd of matrix operator}
		Let $v$ be an $s$-moderate weight, and let $1 \leq p_{i}, q_{i}, s_{t}<\infty$, for $1 \leq i \leq r,$ and $t \in \{1,2\}$, be such that $\frac{1}{p_{1}}+\cdots+\frac{1}{p_{r}}=\frac{1}{s_{1}}$ and $\frac{1}{q_{1}}+\cdots+\frac{1}{q_{r}}=\frac{1}{s_{2}}$. If $\left(a_{i j,m_1 n_1,...,m_r n_r}\right) \in \ell_{\tilde{\Omega}_{s}}^{1}\left(\mathbb{Z}^{2n(r+1)}\right)$, then $\mathcal{O}$ is a bounded operator from $\ell_{\tilde{v}}^{p_{1}, q_{1}}\left(\mathbb{R}^{d}\right) \times \cdots \times \ell_{\tilde{v}}^{p_{r}, q_{r}}\left(\mathbb{R}^{d}\right)$ into $\ell_{\tilde{v}}^{s_{1}, s_{2}}\left(\mathbb{R}^{d}\right)$. In particular, if $\left(a_{i j,m_1 n_1,...,m_r n_r}\right) \in \ell^{1}\left(\mathbb{Z}^{2d(r+1)}\right),$ then $\mathcal{O}$ is a bounded operator from $\ell_{v}^{p_{1}, q_{1}}\left(\mathbb{R}^{d}\right) \times \cdots \times \ell_{v}^{p_{r}, q_{r}}\left(\mathbb{R}^{d}\right)$ into $\ell_{v}^{s_{1}, s_{2}}\left(\mathbb{R}^{d}\right)$.
	\end{theorem}
	The following   result shows that a multilinear integral operator with kernel in the modulation space $M_{\Omega_{s}}^{1}$ gives rise to a bounded operator.
	\begin{theorem}\label{bdd of kernel operator}
		Let $v$ be an $s$-moderate weight, and let $1 \leq p_{i}, q_{i}, s_{t}<\infty,$ for $1 \leq i \leq r,$ and $t \in \{1,2\}$, be such that $\frac{1}{p_{1}}+\cdots+\frac{1}{p_{r}}=\frac{1}{s_{1}}$ and $\frac{1}{q_{1}}+\cdots+\frac{1}{q_{r}}=\frac{1}{s_{2}}$. If $K \in M_{\Omega_{s}}^{1}\left(\mathbb{R}^{d(r+1)}\right)$, then the multilinear integral operator $B_{K}$ defined by \eqref{definition of B_{K}} can be extended as a bounded operator from $M_{v}^{p_{1}, q_{1}}\left(\mathbb{R}^{d}\right) \times \cdots \times M_{v}^{p_{r}, q_{r}}\left(\mathbb{R}^{d}\right)$ into $M_{v}^{s_{1}, s_{2}}\left(\mathbb{R}^{d}\right)$.
	\end{theorem}
	\begin{proof}
		Let $f_1,f_2,...,f_r,h \in \mathcal{S} (\mathbb{R}^d)$. Then,  using the  Gabor series expansion, we have 
		$$
		\begin{aligned}
			&f_{k}=\sum_{m_{k}, n_{k}}\left\langle f_{k}, M_{\beta n_{k}} T_{\alpha m_{k}} \phi\right\rangle M_{\beta n_{k}} T_{\alpha m_{k}} \gamma, \quad 1 \leq k \leq r,\\
			&h=\sum_{i, j}\left\langle h, M_{\beta j} T_{\alpha i} \phi\right\rangle M_{\beta j} T_{\alpha i} \gamma,
		\end{aligned}
		$$
		where $\phi$ and $\gamma$ are dual Gabor frames. By \cite[Proposition 1]{Okoudjou}, the matrix defined by \eqref{infinite matrix form} belongs to $l^{1}_{\tilde{\Omega}_{s}}$, since $K \in M_{\Omega_{s}}^{1}$. Therefore, by Theorem \ref{bdd of matrix operator}, we have the following estimates:
		$$
		\begin{aligned}
			&\left|\langle B_{K}(f_1,f_2,...,f_r ),h\rangle\right|\\
			&=\left|\sum_{i, j} \sum_{m_{1}, n_{1}} \cdots \sum_{m_{r}, n_{r}} b_{i j,m_1 n_1,...,m_r n_r} \overline{\left\langle h, M_{\beta j} T_{\alpha i} \phi\right\rangle} \left\langle f_{1}, M_{\beta n_{1}} T_{\alpha m_{1}} \phi\right\rangle \cdots  \left\langle f_{r}, M_{\beta n_{r}} T_{\alpha m_{r}} \phi\right\rangle\right| \\
			& \leq C\left\|b_{i j,m_1 n_1,...,m_r n_r} \right\|_{l^{1}_{\tilde{\Omega}_{s}}}\left\|\left\langle f_{1}, M_{\beta n_{1}} T_{\alpha m_{1}} \phi\right\rangle \right\|_{l_{\tilde{v}}^{p_{1}, q_{1}}} \cdots \left\| \left\langle f_{r}, M_{\beta n_{r}} T_{\alpha m_{r}} \phi\right\rangle\right\|_{l_{\tilde{v}}^{p_{r}, q_{r}}}\left\|\left\langle h, M_{\beta j} T_{\alpha i} \phi\right\rangle\right\|_{l_{1 / \tilde{v}}^{s_{1}^{\prime}, s_{2}^{\prime}}} \\
			& \leq C\|K\|_{M_{\Omega_{s}}^{1}}\|f_{1}\|_{M_{v}^{p_{1}, q_{1}}}\cdots \|f_{r}\|_{M_{v}^{p_{r}, q_{r}}} \|h\|_{M_{1 / v}^{s_{1}^{\prime}, s_{2}^{\prime}}},
		\end{aligned}
		$$
		where $s_{1}^{\prime}, s_{2}^{\prime}$ are the dual indices of $s_1, s_2,$ respectively. Thus, by the duality, we obtain
		$$
		\left\|B_{K}(f_1,f_2,...,f_r )\right\|_{M_{v}^{s_{1}, s_{2}}} \leq C\|K\|_{M_{\Omega s}^{1}}\|f_{1}\|_{M_{v}^{p_{1}, q_{1}}}\cdots \|f_{1}\|_{M_{v}^{p_{r}, q_{r}}}.
		$$
		Then by a  standard density arguments and   the fact that $\mathbb{S}\left(\mathbb{R}^{d}\right)$ is dense in $M_{v}^{p, q}$ for $1 \leq p, q<\infty$, we get our desire  result.
	\end{proof}
	An immediate application of the above result together with Proposition \ref{kernel expension} and  Theorem \ref{relation of kernel and symbol}  yield  a sufficient condition on the symbol so that the corresponding Fourier integral operator  is bounded on products of modulation spaces.
	\begin{theorem}\label{bdd of symbol operator}
		Let $v$ be an $s$-moderate weight, and let $1 \leq p_{i}, q_{i}, s_{t}<\infty,$ for $1 \leq i \leq r,$ and $t \in \{1,2\}$, be  such that $\frac{1}{p_{1}}+\cdots+\frac{1}{p_{r}}=\frac{1}{s_{1}}$ and $\frac{1}{q_{1}}+\cdots+\frac{1}{q_{r}}=\frac{1}{s_{2}}$.  Let $\sigma_{0}(t_1, t_2,...,t_r, t_{r+1})=\sigma(t_1, t_2,...,t_r, t_{r+1}) $ $\times e^{2\pi i \Phi_{1}(t_1,t_2)}\cdots e^{2\pi i \Phi_{r}(t_1, t_{r+1})}$. If $\sigma_{0} \in M_{\Omega_{s}^{B}}^{1}\left(\mathbb{R}^{d(r+1)}\right)$,  then the corresponding  Fourier integral operator $T_{\sigma_{0}}$  can be extended as a bounded operator from $M_{v}^{p_{1}, q_{1}}\left(\mathbb{R}^{d}\right) \times \cdots \times M_{v}^{p_{r}, q_{r}}\left(\mathbb{R}^{d}\right)$ into $M_{v}^{s_{1}, s_{2}}\left(\mathbb{R}^{d}\right)$.
	\end{theorem}
	\begin{proof}
		From Theorem \ref{relation of kernel and symbol}, we know that  $\sigma_{0} \in M_{\Omega_{s}^{B}}^{1}$ if and only if $K \in M_{\Omega_{s}}^{1},$ where $K$ is the kernel of the corresponding integral operator. Therefore  the result follows from Theorem \ref{bdd of kernel operator}.
	\end{proof}
	Further, if we assume that $v=\omega_{0} \equiv 1$, and that $p_{1}=q_{1}=a_{1}, \cdots p_{r}=q_{r}=a_{r} ~~(\text{hence } s_1=s_2=s)$, then we obtain the following result immediately.
	\begin{corollary}
		Let $2 \leq a_{i}<\infty,$ for $1 \leq i \leq r$, and $s \in \{1,2\}$, be such that $\frac{1}{a_{1}}+\cdots+\frac{1}{a_{r}}=\frac{1}{s}$. Let  $\sigma_{0} \in M^{1}\left(\mathbb{R}^{d(r+1)}\right)$, then $T_{\sigma_{0}}$ can be extended to a bounded operator from $L^{a_1}\left(\mathbb{R}^{d}\right) \times \cdots \times L^{a_{r}}\left(\mathbb{R}^{d}\right)$ into $L^{s}\left(\mathbb{R}^{d}\right)$. In particular, if $\sigma_{0} \in M^{1}\left(\mathbb{R}^{d(r+1)}\right)$, then $T_{\sigma}$ has a bounded extension from $\underbrace{L^{2}\left(\mathbb{R}^{d}\right) \times \cdots \times L^{2}\left(\mathbb{R}^{d}\right)}_{r-times}$ into $L^{\frac{2}{r}}\left(\mathbb{R}^{d}\right.)$.
	\end{corollary}
	\begin{proof}
		Notice that 	for the given range of $a_1,...,a_{r},$ we have the   continuous embeddings; $L^{a_{i}} \subset M^{a_{i}}, $ for $1 \leq i \leq r,$ and so $L^{a_1}\left(\mathbb{R}^{d}\right) \times \cdots \times L^{a_{r}}\left(\mathbb{R}^{d}\right) \subset M^{a_1} \times \cdots \times M^{a_r}$. Moreover, since $1 \leq s \leq 2$, we have that $M^{s} \subset L^{s}$ (see \cite{GK}). Using these continuous embeddings and Theorem \ref{bdd of symbol operator}, we have our  desired result.
	\end{proof}
	\section{Multilinear Fourier integral operators on $\mathbb{T}^d$}\label{sec5}
	In this section we study multilinear   Fourier integral operators on the torus  $\mathbb{T}^d$. First, we will define the following symbol class which is just the multilinear version periodic  H\"ormander symbol class on $\mathbb{T}^d.$
	\begin{defn}\label{symbol class on T^n}
		Let $0 \leq \rho, \delta \leq 1$. We say a symbol $S_{\rho,\delta}^{m_1,m_2,...,m_r}\left(\mathbb{T}^d \times\underbrace{\mathbb{Z}^d \times  \mathbb{Z}^d \times \cdots \times \mathbb{Z}^d}_{r-times}\right),$ where $m_1,m_2,...,m_r \in \mathbb{R},$ if for all multi-indices $\alpha, \beta_1, \beta_2,...,\beta_r$,  there exists a positive constant $C_{\alpha, \beta_1, \beta_2,...,\beta_r}$ such that
		$$|(\partial_{x}^\alpha \triangle_{\xi_{1}}^{\beta_{1}} \cdots \triangle_{\xi_{r}}^{\beta_{r}}\sigma)(x, \xi_1, \xi_2,...,\xi_r) | \leq C_{\beta_1, \beta_2,...,\beta_r} \prod_{k+1}^{r}\left(\langle \xi_{k} \rangle^{m_{k}-\rho|\beta_k|+\alpha |\delta|}\right) ,\quad  x,\xi_1, \xi_2,...,\xi_r \in \mathbb{R}^d.$$
	\end{defn}
	The corresponding periodic multilinear Fourier integral operator  for  a  real-valued phase functions $\Phi_{i} : \mathbb{T}^d \times \mathbb{Z}^{d} \rightarrow \mathbb{R},$   such that   $\Phi_{i}$'s are linear in the second variable for $1 \leq i \leq r$  is defined  by
	\begin{equation}\label{definition of periodic multilinear}
		T_{\sigma}(f )(x)=\sum_{k\in  \mathbb{Z}^{dr}} e^{2\pi i [\Phi_{1}(x, k_1)+\cdots+\Phi_{r}(x, k_r)]}  \sigma(x, k) \hat{f_1}(k_1) \cdots \hat{f_r}(k_r), \quad x \in \mathbb{T}^d,
	\end{equation}
	where   $ k=(k_1,...,k_r)\in \mathbb{Z}^{dr}, f=( f_1,f_2,...,f_r) \in  C^{\infty} (\mathbb{T}^d)^r$, and $$\hat{f_i}(k_i) = \int_{\mathbb{T}^d} e^{-2 \pi i \eta\cdot k_{i}} f_{i}(\eta), d \eta, \quad k_i\in \mathbb{Z}^d,$$  is the periodic  Fourier transform of $f_i$.
	We start with the following well defined result. 
	\begin{prop}\label{periodic operator well defined}
		Let $\sigma(x, \xi_1, \xi_2,...,\xi_r) \in S^{m_1,m_2,...,m_r}\big(\mathbb{T}^d \times\underbrace{\mathbb{Z}^d \times  \mathbb{Z}^d \times \cdots \times \mathbb{Z}^d}_{r-times}\big),$ where $m_1,m_2,...,$ $m_r  \in \mathbb{R}$, and assume that there exist positive constants $C_1,C_2,...,C_r$ satisfies the  conditions $\left|\partial_{x}^{\alpha} e^{i \Phi_{l}(x,\xi_l)}\right| \leq C_l$, for every multi-index $\alpha$ and $1 \leq l \leq r$. 	Then   the corresponding periodic multilinear Fourier integral operator   $T_{\sigma}$ is well defined.
	\end{prop}
	\begin{proof}
		For  $f=( f_1,f_2,...,f_r) \in  C^{\infty} (\mathbb{T}^d)^r$, the periodic multilinear Fourier integral operator defined as 
		\begin{equation*} 
			T_{\sigma}(f )(x)=\sum_{\xi \in  \mathbb{Z}^{dr}} e^{2\pi i [\Phi_{1}(x, \xi_1)+\cdots+\Phi_{r}(x, \xi_r)]}  \sigma(x, \xi) \hat{f_1}(\xi_1) \cdots \hat{f_r}(\xi_r), \quad x \in \mathbb{T}^d,
		\end{equation*}
		where   $ \xi=(\xi_1,...,\xi_r)\in \mathbb{Z}^{dr}$. This can also be written as 
		$$
		T_{\sigma}(f )(x)=\sum_{\xi_1 \in  \mathbb{Z}^{d}} e^{2\pi i \Phi_{1}(x, \xi_1)}  \sigma_1(x,\xi_1) \hat{f_1}(\xi_1),
		$$
		where 
		\begin{align}\label{symbol for T}
			\sigma_1(x,\xi_1)=\sum_{\xi_2,...,\xi_r \in  \mathbb{Z}^{d}} e^{2\pi i [\Phi_{2}(x, \xi_2)+\cdots+\Phi_{r}(x, \xi_r)]}  \sigma(x, \xi) \hat{f_2}(\xi_2) \cdots \hat{f_r}(\xi_r).
		\end{align}
		For a moment, we assume that $\sigma_1(x,\xi_1) \in S^{m_1}\left(\mathbb{T}^d \times \mathbb{Z}^d\right).$ Then  using the fact that $\Phi(x,\xi_1)$ is linear in $\xi_1$ and  the definition of forward difference $\triangle$, we obtained
		$$\left|T_{\sigma}(f)(x) \right| \leq \sum_{\xi_1 \in  \mathbb{Z}^{d}} B_{1} \langle \xi\rangle^{m_1-k_1\rho},$$ 
		where $ B_1 > 0.$
		Now, choose $k_1$ sufficiently large such that $m_1-k_1\rho << -d.$ This gives us that the above sum is absolutely convergent. Thus it remains to prove that $\sigma_1(x,\xi_1) \in S^{m_1}\left(\mathbb{T}^d \times \mathbb{Z}^d\right).$  From (\ref{symbol for T}), again we can write 
		$$
		\begin{aligned}
			\sigma_1(x,\xi_1)&=\sum_{\xi_2,...,\xi_r \in  \mathbb{Z}^{d}} e^{2\pi i [\Phi_{2}(x, \xi_2)+\cdots+\Phi_{r}(x, \xi_r)]}  \sigma(x, \xi) \hat{f_2}(\xi_2) \cdots \hat{f_r}(\xi_r)\\
			&=\sum_{\xi_2 \in  \mathbb{Z}^{d}} e^{2\pi i \Phi_{2}(x, \xi_2)}  \sigma_2(x,\xi_1,\xi_2) \hat{f_2}(\xi_2) ,\\
		\end{aligned}
		$$
		where 
		$$
		\sigma_2(x,\xi_1,\xi_2)=\sum_{\xi_3,...,\xi_r \in  \mathbb{Z}^{d}} e^{2\pi i [\Phi_{3}(x, \xi_3)+\cdots+\Phi_{r}(x, \xi_r)]}  \sigma(x, \xi ) \hat{f_3}(\xi_3) \cdots \hat{f_r}(\xi_r).
		$$
		Similarly, as above, for the  moment, if we assume that $\sigma_2(x,\xi_1,\xi_2) \in S^{m_1,m_2}\left(\mathbb{T}^d \times \mathbb{Z}^d \times \mathbb{Z}^d\right).$ Again,   using the fact that $\Phi(x,\xi_2)$ is linear in $\xi_2$ and  the definition of forward difference $\triangle$ along with the     given condition, we obtained $\sigma_1(x,\xi_1) \in S^{m_1}\left(\mathbb{T}^d \times \mathbb{Z}^d\right).$ Now, to prove $\sigma_2(x,\xi_1,\xi_2) \in S^{m_1,m_2}\left(\mathbb{T}^d \times \mathbb{Z}^d \times \mathbb{Z}^d\right),$ we again use the similar techniques developed as above. In the end, we obtained $\sigma_{r}(x, \xi_1, \xi_2,...,\xi_{r-1}) \in S^{m_1,m_2,...,m_r}\big(\mathbb{T}^d \times\underbrace{\mathbb{Z}^d \times  \mathbb{Z}^d \times \cdots \times \mathbb{Z}^d}_{(r-1)-times}\big)$, where
		$$
		\sigma_{r}(x, \xi_1, \xi_2,...,\xi_{r-1}) = \sum_{\xi_r \in \mathbb{Z}^{d}} e^{2\pi i \Phi_{r}(x, \xi_r)} \sigma(x, \xi) \hat{f_r}(\xi_r) ,
		$$
		and this completes the proof.
	\end{proof}
	Now we define periodic multilinear integral operator associated with kernel $K \in \mathcal{D}^{\prime}\left(\mathbb{T}^{n(r+1)}\right)$, which is a mapping $B_{K}$ from $\underbrace{C^{\infty}\left(\mathbb{T}^{d}\right) \times \cdots \times  C^{\infty}\left(\mathbb{T}^{d}\right)}_{r-times}$ to $\mathcal{D}^{\prime}\left(\mathbb{T}^{d}\right)$ by
	\begin{equation}\label{periodic definition of B_{K}}
		B_{K}(f_1,f_2,...,f_r )(x)=\int_{\mathbb{T}^{dr}} K(x, y_1,...,y_r) f_1(y_1) \cdots f_r(y_r) d y_1 \cdots d y_r,
	\end{equation}
	where $f_1,f_2,...,f_r \in C^{\infty}(\mathbb{T}^d)$. In the next result we establish the relationship between a periodic multilinear integral operator
	and a periodic multilinear Fourier integral operator defined by \eqref{definition of periodic multilinear}.
	\begin{theorem}\label{periodic kernel expension}
		Let $T_{\sigma}$ be a periodic multilinear Fourier integral operator associated with  symbol $\sigma$. Then $T_{\sigma}$ coincides with a periodic multilinear integral operator $B_{K}$ with kernel
		$$K(x, \xi)=\mathcal{F}_{r+1} \mathcal{F}_{r} \cdots \mathcal{F}_{3} \mathcal{F}_{2}(\sigma(x, \xi  ) e^{2\pi i \Phi_{1}(x, \xi_1)} \cdots e^{2\pi i \Phi_{r}(x, \xi_r)}).$$
		where $\xi=(\xi_1, \xi_2,...,\xi_r)$ and  $ \mathcal{F}_{j}$ denotes the Fourier transform in the $j^{th}$ variable.
	\end{theorem}
	\begin{proof}
		For $f=( f_1,f_2,...,f_r)$ with $f_i  \in C^{\infty}(\mathbb{T}^d)$, we have 
		$$
		\begin{aligned}
			T_{\sigma}(f )(x)
			&=\sum_{\xi_1,...,\xi_r \in  \mathbb{Z}^{d}} e^{2\pi i [\Phi_{1}(x, \xi_1)+\cdots+\Phi_{r}(x, \xi_r)]}  \sigma(x, \xi_1,...,\xi_r) \hat{f_1}(\xi_1) \cdots \hat{f_r}(\xi_r) \\
			&=\sum_{\xi_1,...,\xi_r \in  \mathbb{Z}^{d}} \int_{\mathbb{T}^{dr}} e^{2\pi i [\Phi_{1}(x, \xi_1)-\xi_1 \cdot y_1+\cdots+\Phi_{r}(x, \xi_r)-\xi_r \cdot y_r]}  \sigma(x, \xi ) f(y_1) \cdots f(y_r) \;d y_1 \cdots d y_r\\
			&=\int_{\mathbb{T}^{dr}} K(x, y_1,\cdots,y_r) f_1(y_1) \cdots f_r(y_r) \;d y_1 \cdots d y_r\\
			&=B_{K}(f )(x),
		\end{aligned}
		$$
		where
		$$
		\begin{aligned}
			K(x, y_1,\cdots,y_r) &= \sum_{\xi_1,...,\xi_r \in  \mathbb{Z}^{d}} e^{2\pi i [\Phi_{1}(x, \xi_1)-\xi_1 \cdot y_1+\cdots+\Phi_{r}(x, \xi_r)-\xi_r \cdot y_r]}  \sigma(x, \xi )\\
			&=\mathcal{F}_{r+1} \mathcal{F}_{r} \cdots \mathcal{F}_{3} \mathcal{F}_{2}\left(\sigma(x, y_1, y_2,...,y_r ) e^{2\pi i \Phi_{1}(x, y_1)} \cdots e^{2\pi i \Phi_{r}(x, y_r)}\right).\\
		\end{aligned}
		$$
	\end{proof}
	For $\xi=(\xi_1, \xi_2,...,\xi_r)$,   suppose that $\sigma(x, \xi ) \in S^{m_1,m_2,...,m_r}\big(\mathbb{T}^d \times\underbrace{\mathbb{Z}^d \times  \mathbb{Z}^d \times \cdots \times \mathbb{Z}^d}_{r-times}\big)$ be such that $\sigma(x, \xi_1, \xi_2,...,\xi_r)  e^{2\pi i \Phi_{1}(x, \xi_1)} \cdots e^{2\pi i \Phi_{r}(x, \xi_r)}$ is invariant under fourier transform on Modulation space $M^{1}\left(\mathbb{T}^{d(r+1)}\right).$ Then, we obtained the following result.
	\begin{corollary}\label{relation of kernel and symbol in periodiccase}
		Let $T_{\sigma}$ be a periodic multilinear Fourier integral operator associated with  symbol $\sigma$. Then 	$\sigma(x, y_1, y_2,...,y_r) e^{2\pi i \Phi_{1}(x, y_1)} \cdots e^{2\pi i \Phi_{r}(x, y_r)} \in M^{1}\left(\mathbb{T}^{n(r+1)}\right)$ if and only if $K(x, y_1, y_2,...,y_r)=\mathcal{F}_{r+1} \mathcal{F}_{r} \cdots \mathcal{F}_{2}\left(\sigma(x, y_1, y_2,...,y_r ) e^{2\pi i \Phi_{1}(x, y_1)} \cdots e^{2\pi i \Phi_{r}(x, y_r)}\right)\in M^{1}\left(\mathbb{T}^{d(r+1)}\right).$
	\end{corollary}
	\section{Continuity of bilinear pseudo-differential operator }\label{sec6}
	In this section, we study the boundedness of bilinear pseudo-differential operators on modulation spaces for certain symbol classes, namely $\mathbf{SG}$-class $\mathbf{SG}^{m_{1}, m_{2}, m_{3}}$. A symbol $\sigma$ on $\mathbb{R}^{3d}$ is in $\mathbf{SG}^{m_1,m_2,m_3}$, where $m_1,m_2,m_3 \in \mathbb{R},$ if for all multi-indices $\alpha,\beta$ and $\gamma,$ there exists a positive constant $C_{\alpha,\beta,\gamma}$ such that it satisfies the following condition:
	$$
		\left|\partial_{x}^{\alpha} \partial_{\xi}^{\beta} \partial_{\eta}^{\gamma} \sigma(x, \xi, \eta)\right| \leq C_{\alpha, \beta, \gamma}\langle x \rangle^{m_{3}-|\alpha|}\langle \xi\rangle^{m_{1}-|\beta|} \langle \eta\rangle^{m_{2}-|\gamma|}, \quad x,\xi,\eta \in \mathbb{R}^d.
	$$
	Corresponding to the symbol $\sigma$, a bilinear pseudo-differential operator $T_\sigma: \mathcal{S}\left(\mathbb{R}^d\right) \times \mathcal{S}\left(\mathbb{R}^d\right)\to \mathcal{S}^{\prime}\left(\mathbb{R}^d\right) $  is defined as 
	$$
	T_\sigma(f, g)(x)=\int_{\mathbb{R}^d} \int_{\mathbb{R}^d} e^{2 \pi i x \cdot(\xi+\eta)} \sigma(x, \xi, \eta) \hat{f}(\xi) \hat{g}(\eta) \;  d \xi d \eta, \qquad x\in \mathbb{R}^d,
	$$
	for $f, g \in \mathcal{S}\left(\mathbb{R}^d\right)$. 
	
	For $s_{1}, s_{2} \in \mathbb{R}$, we define $v_{s_{1}, s_{2}}(x,y)=\langle x \rangle^{s_{2}}\langle y \rangle^{s_{1}}\langle\eta\rangle^{s_{2}}$. 
	From now onwards, for  given $N_1, N_2,$ and $ N_3\in \mathbb{N}$,   we consider rough symbols $\sigma$   on $\mathbb{R}^{3d}$ satisfying estimates of the type
	\begin{equation}\label{symbol estimate for continuity result}
		\left|\partial_{x}^{\alpha} \partial_{\xi}^{\beta} \partial_{\eta}^{\gamma} \sigma(x, \xi, \eta)\right| \leq C_{\alpha, \beta, \gamma}\langle x \rangle^{m_{3}}\langle \xi\rangle^{m_{1}} \langle \eta\rangle^{m_{2}}, \quad|\alpha| \leq 2 N_{3},|\beta| \leq 2 N_{1},|\gamma| \leq 2 N_{2},
	\end{equation}
with $\partial_{x}^{\alpha} \partial_{\xi}^{\beta} \partial_{\eta}^{\gamma}$ standing for distributional derivatives. Then our main result of this section is as follows.
	\begin{theorem}\label{bdd of symbol in modulatiobn space}
		For $s_{1}<<0, s_{2}>0$, let $\mu \in \mathcal{M}_{v_{s_{1}, s_{2}}}$.
		Consider a symbol $\sigma$ satisfying \eqref{symbol estimate for continuity result}, with $N_{1}>\frac{s_2 +d}{2}$, $N_{2}>\frac{d}{2}$, and  $\frac{|m_2|+d}{2} < N_{3} < \frac{-s_1-|m_1|-d}{2}$. Then, for every $1 \leq p, q \leq \infty,$ the corresponding pseudo-differential operatr $T_\sigma$ extends to a continuous operator from $\tilde{M}_{\mu}^{p, q} \otimes \tilde{M}_{\mu}^{p, q}$ to $\tilde{M}_{\mu v_{-m_{1}-m_{2}-2 N_{3},-m_3}}^{p, q}$.
	\end{theorem}
	To prove the above theorem, we require some technical preparation. First, we prove an almost diagonalization result for bilinear Fourier integral operators in the case of regular symbols with respect to a Gabor frame, and then obtained the result for bilinear pseudo-differential operators as a special case. Assume that the function $\Phi(x,\xi, \eta )$ fulfill the following properties:
	\begin{enumerate}[(i)]
		\item \label{1}$\Phi \in C^{\infty}(\mathbb{R}^{3d}),$
		\item \label{2} for $z=(x, \xi, \eta)$, 
		$ | \partial^{\alpha} \Phi(z)| \leq c_{\alpha}, ~ |\alpha| \geq 2,$
		\item there exist $\delta \geq 0$ such that $|\det \partial^3_{x,\xi, \eta} \Phi(x, \xi, \eta)| \geq \delta$.
	\end{enumerate}
	Here, for a given $N \in \mathbb{N}$, we consider $\sigma$ on $\mathbb{R}^{3n}$ satisfying
	\begin{equation}\label{symbol derivative estimate}
		|\partial^{\alpha}_z \sigma(z)| \leq c_{\alpha}
		\text{ a.e } z=(x, \xi, \eta) \in \mathbb{R}^{3d},\quad |\alpha| \leq 2N, \end{equation} 
	where $\partial^{\alpha}_{z}$ denotes distributional derivatives.
	The following result is the analogus of \cite[Theorem 5.1]{EC&FN&LR} and can be proved in similar lines with the help of bilinear interpolation theory which can be found in \cite{bergh}. 
	\begin{prop}\label{kernel bdd in bilinear case}
		Consider an operator 
		$$(K(c,e))_{m^{\prime}, n^{\prime}} = \sum_{m,n,m_0,n_0} K_{m^{\prime},n^{\prime}, m, n, m_{0}, n_{0}} c_{m, n} e_{m_{0}, n_{0}},$$
		which is defined on sequences on the lattice
		$\Lambda^{\prime} = \left(\alpha \mathbb{Z}^d \times \beta \mathbb{Z}^d\right) \times \left(\alpha \mathbb{Z}^d \times \beta \mathbb{Z}^d\right).$
		\begin{enumerate}[(i)]
			\item  If $K \in \ell^{\infty}_n \ell^{\infty}_{n_0} \ell^1_{n^{\prime}} \ell^{\infty}_{m^{\prime}} \ell^1_{m} \ell^1_{m_0}$ then $K: \ell^1_n \ell^{\infty}_m \otimes \ell^1_{n_{0}} \ell^{\infty}_{m_{0}} \rightarrow \ell^{1}_{n^{\prime}} \ell^{\infty}_{m^{\prime}}$ is a continuous operator. 
			\item  If $K \in \ell^{\infty}_{n^{\prime}} \ell^1_{n_{0}} \ell^{1}_{n} \ell^{\infty}_{m} \ell^{\infty}_{m_0} \ell^1_{m^{'}}$ then $K: \ell^{\infty}_n \ell^{1}_m \otimes \ell^{\infty}_{n_{0}} \ell^{1}_{m_{0}} \rightarrow \ell^{\infty}_{n^{\prime}} \ell^{1}_{m^{\prime}}$ is a continuous operator. 
			\item If $K \in \ell^{\infty}_n \ell^{\infty}_{n_0} \ell^1_{n^{\prime}} \ell^{\infty}_{m^{\prime}} \ell^1_{m} \ell^1_{m_0} \cap   \ell^{\infty}_{n^{\prime}} \ell^{1}_{n} \ell^1_{n_{0}} \ell^{\infty}_{m} \ell^{\infty}_{m_0} \ell^1_{m^{\prime}}$ and   $K \in \ell^{\infty}_{n^{\prime}} \ell^{\infty}_{m^{\prime}} \ell^{1}_{n} \ell^{1}_{m} \ell^1_{n_{0}}  \ell^1_{m_{0}} \cap \ell^{\infty}_{n} \ell^{\infty}_{m} \ell^{\infty}_{n_0} \ell^{\infty}_{m_0} \ell^{1}_{n^{\prime}} \ell^{1}_{m^{\prime}},$ then $K: {\ell}^{p,q} \otimes  {\ell}^{p,q} \rightarrow  {\ell}^{p,q}$ is a continuous operator for all $1 \leq p,q \leq \infty$, where $ {\ell}^{p,q}={\ell}^{q}_{n} {\ell}^{p}_{m}.$
			\item Assume the hypothesis in (iii). Then $K: \tilde{\ell}^{p,q} \otimes \tilde{\ell}^{p,q} \rightarrow \tilde{\ell}^{p,q}$ is a continuous operator for all $1 \leq p,q \leq \infty$. 
		\end{enumerate}
	\end{prop}
	Now, we give the decay properties of the matrix of the bilinear Fourier integral operator $T$ with respect to a Gabor frame which plays an important role to obtain the main result of this section.
	

	\begin{theorem}\label{estimate for bilinear FIO}
		Consider a phase function satisfying (\ref{1}) and (\ref{2}) and a symbol satisfying \eqref{symbol derivative estimate}. Then there exists  a positive constant $C_N $ such that 
		\begin{equation}\label{estimate of symbol in bilinear FIO case}
			|\langle T_{\sigma}(g_{m,n},g_{m_{0},n_{0}}), g_{m^{\prime},n^{\prime}} \rangle | \leq C_N~ \langle \nabla_z \Phi(m^{\prime},n,n_{0})- (n^{\prime},m,m_{0}) \rangle^{-2N} .
		\end{equation}
	\end{theorem}
	\begin{proof}
		Using the fact that $\left(T_x f\right)^{\wedge}=M_{-x} \hat{f},$   $\left(M_\eta f\right)^{\wedge}=T_\eta \hat{f}$,  and  the commutation relations $T_x M_\eta=e^{-2 \pi i x \eta} M_\eta T_x$,  we can  write
		\begin{align}\label{3}\nonumber
			&	\langle T_{\sigma}(g_{m,n},g_{m_{0},n_{0}}), g_{m^{\prime},n^{\prime}}\rangle\\\nonumber
			&= \int_{\mathbb{R}^{d}}  \int_{\mathbb{R}^{d}}  \int_{\mathbb{R}^{d}}   e^{2\pi i \Phi(x, \xi, \eta)} \sigma(x,\xi, \eta) T_n M_{-m} \hat{g}(\xi) T_{n_{0}} M_{-m_{0}} \hat{g}(\eta) M_{-n^{\prime}} T_{m^{\prime}} \bar{g}(x)\; dx d\xi d\eta\\\nonumber
			&= \int_{\mathbb{R}^{d}}  \int_{\mathbb{R}^{d}}  \int_{\mathbb{R}^{d}}    M_{(0,0,-m_{0})} T_{(0,0,-n_{0})} M_{(0,-m,0)} T_{(0,-n,0)} e^{2\pi i \Phi(x, \xi, \eta)} \sigma(x,\xi, \eta) \\\nonumber
			&\qquad \times M_{-n^{\prime}} T_{m^{\prime}} \bar{g}(x) \hat{g}(\xi) \hat{g}(\eta) \;dx d\xi d\eta\\\nonumber
			&= \int_{\mathbb{R}^{d}}  \int_{\mathbb{R}^{d}}  \int_{\mathbb{R}^{d}}    T_{(-m^{\prime},0,0)} M_{(-n^{\prime},0,0)} T_{(0,0,-n_{0})} M_{(0,-m,0)} T_{(0,-n,0)} e^{2\pi i \Phi(x, \xi, \eta)} \sigma(x,\xi, \eta)
			\\\nonumber
			&\qquad \times M_{-n^{\prime}} T_{m^{\prime}} \bar{g}(x) \hat{g}(\xi) \hat{g}(\eta) \;dx d\xi d\eta\\\nonumber
			&= \int_{\mathbb{R}^{d}}  \int_{\mathbb{R}^{d}}  \int_{\mathbb{R}^{d}}    e^{2\pi i(\Phi(x+m^{\prime}, \xi+n, \eta +n_{0})-(n^{\prime},m,m_{0}).(x+m^{\prime},\xi, \eta)}) \sigma(x+m^{\prime}, \xi+n, \eta+n_{0})\\&\qquad \times \bar{g}(x) \hat{g}(\xi) \hat{g}(\eta) dx d\xi d\eta.
		\end{align}
		Since $\Phi$ is smooth, we expand $\Phi(x, \xi, \eta)$ into a Taylor series around $(m, n, n_{0})$ and obtain
		$$ \Phi(x+m^{\prime}, \xi+n, \eta +n_{0}) = \Phi(m^{\prime}, n,n_{0}) + \nabla_z\Phi(m^{\prime}, n, n_{0}).(x, \xi, \eta) + \Phi_{2,(m^{\prime},n,n_0)} (x,\xi,\eta),$$where reminder is given by $$\Phi_{2,(m^{\prime},n,n_0)} (x,\xi,\eta)= 2 \frac{(x,\xi,\eta)^{\alpha}}{\alpha !} \sum_{|\alpha| = 2} \int_{0}^{1} (1-t) \partial^{\alpha} \Phi((m^{\prime}, n, n_{0})+ t(x, \xi, \eta)) dt . $$
		Also, for any $N \in \mathbb{N}$, we have the following identity
		\begin{equation}\label{estimate of exponential part}
			\begin{aligned}[b]
				&(1-\nabla_{z})^{N} e^{2\pi i \{[\nabla_z \Phi(m^{\prime}, n, n_{0})-(n^{\prime}, m, m_{0})].(x,\xi,\eta)\}}\\
				&= \langle 2 \pi(\nabla_z \Phi(m^{\prime}, n, n_{0})- (n^{\prime}, m, m_{0}) \rangle^{2N} e^{2\pi i \{[\nabla_z \Phi(m^{\prime}, n, n_{0})-(n^{\prime}, m, m_{0})].(x,\xi,\eta)\}} 
			\end{aligned}
		\end{equation}
		Now, using integration by parts and the relation \eqref{estimate of exponential part},  from (\ref{3}), we get 
		$$ 
		\begin{aligned}
			&|\langle T_{\sigma}(g_{m,n},g_{m_{0},n_{0}}), g_{m^{'},n^{'}} \rangle|\\
			&= \left| \int_{\mathbb{R}^{3d}} e^{2\pi i \{[\nabla_z\Phi(m^{\prime}, n, n_{0})-(n^{\prime}, m, m_{0})].(x,\xi,\eta)\}} \sigma(x+m^{\prime}, \xi+n, \eta+n_{0}) e^{2 \pi i \Phi_{2,(m^{\prime},n,n_0)} (x,\xi,\eta)}   \bar{g}(x) \hat{g}(\xi) \hat{g}(\eta) dx d\xi d\eta\right|\\
			&=\frac{1}{\langle 2 \pi(\nabla_z \Phi(m^{\prime}, n, n_{0})- (n^{\prime}, m, m_{0})) \rangle^{2N}} \left|\int_{\mathbb{R}^{3d}} e^{2\pi i \{[\nabla_z \Phi(m^{\prime}, n, n_{0})-(n^{\prime}, m, m_{0})].(x,\xi,\eta)\}}\right. \\
			& \qquad \left.\times (1-\nabla_z)^{N} \left[\sigma(x+m^{\prime}, \xi+n, \eta+n_{0}) e^{2 \pi i \Phi_{2,(m^{\prime},n,n_0)} (x,\xi,\eta)}   \bar{g}(x) \hat{g}(\xi) \hat{g}(\eta) \right]dx d\xi d\eta\right|.
		\end{aligned}
		$$
		By means of the Leibniz's formula the factor
		$$
		(1-\nabla_z)^{N} \left[\sigma(x+m^{\prime}, \xi+n, \eta+n_{0}) e^{2 \pi i \Phi_{2,(m^{\prime},n,n_0)} (x,\xi,\eta)}   \bar{g}(x) \hat{g}(\xi) \hat{g}(\eta) \right]
		$$
		can be expressed as
		$$
		e^{2 \pi i \Phi_{2,\left(m^{\prime}, n,n_{0}\right)}(z)} \sum_{|\alpha|+|\beta|+|\gamma| \leq 2 N} C_{\alpha, \beta \gamma} p\big(\partial^{|\alpha|} \Phi_{2,\left(m^{\prime}, n,n_{0}\right)}\big)(z)  \big(\partial_{z}^{\beta} \sigma\big)\left(z+\left(m^{\prime}, n,n_{0}\right)\right)\cdot \partial_{z}^{\gamma}(\bar{g} \otimes \hat{g}  \otimes \hat{g})(z),
		$$
		where $p\left(\partial^{|\alpha|} \Phi_{2,\left(m^{\prime}, n,n_0\right)}\right)(z)$ is a polynomial made of derivatives of $\Phi_{2,\left(m^{\prime}, n,n_0\right)}$ of order at most $|\alpha|$.
		As a consequence of (\ref{2}), we have $\partial_{z}^{\alpha} \Phi_{2,\left(m^{\prime}, n,n_0\right)}(z)=O\left(\langle z\rangle^{2}\right)$, which combined with the assumption \eqref{symbol derivative estimate} and the hypothesis $g \in \mathcal{S}\left(\mathbb{R}^{d}\right)$ yields the desired estimate.
	\end{proof}
	An immediate application of the above result is the following one.
	\begin{theorem}
		Consider a symbol $\sigma$ satisfying \eqref{symbol estimate for continuity result}. Then there exists $C_{N_{1}, N_{2}, N_{3}}>0$ such that
		\begin{equation}\label{estimate of bilinear pdo}
			\left|\left\langle T_\sigma( g_{m, n},g_{m_{0}, n_{0}}), g_{m^{\prime}, n^{\prime}}\right\rangle\right| \leq C_{N_{1}, N_{2}, N_{3}} \frac{\langle n\rangle^{m_{1}} \langle n_{0}\rangle^{m_{2}} \left\langle m^{\prime}\right\rangle^{m_{3}}}{\left\langle n+n_{0}-n^{\prime}\right\rangle^{2 N_{3}} \left\langle m-m^{\prime}\right\rangle^{2 N_{1}} \left\langle m_{0}-m^{\prime}\right\rangle^{2 N_{2}}} .
		\end{equation}
	\end{theorem}
	\begin{proof}
		The proof is essentially a particular case of Theorem \ref{estimate for bilinear FIO}, so we only give an outline of the main ideas. An explicit computation shows that
		$$
		\begin{aligned}
			& \left|\left\langle T_\sigma ( g_{m, n},g_{m_{0}, n_{0}}), g_{m^{\prime}, n^{\prime}}\right\rangle\right| \\
			& \quad=\left|\int_{\mathbb{R}^{3d}} e^{2 \pi i\left[ x\left(n+n_{0}-n^{\prime}\right)-\xi(m-m^{\prime})-\eta\left(m_{0}-m^{\prime}\right)\right]}\left[e^{2 \pi i x \cdot (\xi+\eta)} \sigma\left(x+m^{\prime},\xi+n, \eta+n_{0}\right)\right] \bar{g}(x)  \hat{g}(\xi) \hat{g}(\eta) d x d \xi d \eta\right| .
		\end{aligned}
		$$
		Then one uses the identity
		$$
		\begin{aligned}
			& \left(1-\Delta_{x}\right)^{N_{3}}\left(1-\Delta_{\xi}\right)^{N_{1}}  \left(1-\Delta_{\eta}\right)^{N_{2}} e^{2 \pi i\left[ x\left(n+n_{0}-n^{\prime}\right)-\xi(m-m^{\prime})-\eta\left(m_{0}-m^{\prime}\right)\right]} \\
			& =\left\langle 2 \pi\left(n+n_{0}-n^{\prime}\right)\right\rangle^{2 N_{3}}\left\langle 2 \pi\left(m-m^{\prime}\right)\right\rangle^{2 N_{1}} \left\langle 2 \pi\left(m_{0}-m^{\prime}\right)\right\rangle^{2 N_{2}} e^{2 \pi i\left[ x\left(n+n_{0}-n^{\prime}\right)-\xi(m-m^{\prime})-\eta\left(m_{0}-m^{\prime}\right)\right]},
		\end{aligned}
		$$
		and integration by parts. Since $g \in \mathcal{S}$, the estimates \eqref{symbol estimate for continuity result} combined with Peetre's inequality $\langle z+w\rangle^{s} \leq\langle z\rangle^{s}\langle w\rangle^{|s|}$ gives us the required estimate \eqref{estimate of bilinear pdo}.
	\end{proof}
	Now we are in a position to prove our main result of this section. 
	\begin{proof}[Proof of Theorem \ref{bdd of symbol in modulatiobn space}]
		Consider a normalized tight bilinear frame $\mathcal{G}(g \otimes g, \alpha, \beta)$ with $g \in \mathcal{S}\left(\mathbb{R}^{d}\right)$. From Theorem \ref{gabor frame bdd},   in order to show the boundedness of $T_\sigma $ from $\tilde{M}_{\mu}^{p, q} \otimes \tilde{M}_{\mu}^{p, q}$ to $\tilde{M}_{\mu v_{-m_{1}-m_{2}-2 N_{3},-m_3}}^{p, q}$, it is sufficient to prove the boundedness of the infinite matrix
		$$
		K_{m^{\prime}, n^{\prime}, m, n,m_0,n_0}=\left\langle T_\sigma (g_{m, n},g_{m_{0}, n_{0}}), g_{m^{\prime}, n^{\prime}}\right\rangle \frac{\mu\left(m^{\prime}, n^{\prime}\right)}{\left\langle m^{\prime}\right\rangle^{m_{3}}\left\langle n^{\prime}\right\rangle^{m_{1}} \left\langle n^{\prime}\right\rangle^{m_{2}}\left\langle n^{\prime}\right\rangle^{2N_{3}}\mu(m, n)}
		$$
		from $\tilde{\ell}^{p, q} \otimes \tilde{\ell}^{p, q}$ into $\tilde{\ell}^{p, q}$. The estimate \eqref{estimate of bilinear pdo} and the assumption $\mu \in \mathcal{M}_{v_{s_{1}, s_{2}}}$ combined with Petree's inequality yield
		$$
		\left|K_{m^{\prime}, n^{\prime}, m, n,m_0,n_0}\right| \lesssim\left\langle n-n^{\prime}\right\rangle^{s_{1}+\left|m_{1}\right|+2 N_{3}} \left\langle n_{0}-n^{\prime}\right\rangle^{\left|m_{2}\right|-2 N_{3}} \left\langle m-m^{\prime}\right\rangle^{s_{2}-2 N_{1}} \left\langle m_{0}-m^{\prime}\right\rangle^{-2 N_{2}}.
		$$
		Now, because of the choice of $N_{1}, N_{2}, N_{3}$, one can deduce that $K \in \ell^{\infty}_n \ell^{\infty}_{n_0} \ell^1_{n^{\prime}} \ell^{\infty}_{m^{\prime}} \ell^1_{m} \ell^1_{m_0} \cap   \ell^{\infty}_{n^{\prime}} \ell^{1}_{n} \ell^1_{n_{0}} \ell^{\infty}_{m} \ell^{\infty}_{m_0} \ell^1_{m^{\prime}},$ and $K \in \ell^{\infty}_{n^{\prime}} \ell^{\infty}_{m^{\prime}} \ell^{1}_{n} \ell^{1}_{m} \ell^1_{n_{0}}  \ell^1_{m_{0}} \cap \ell^{\infty}_{n} \ell^{\infty}_{m} \ell^{\infty}_{n_0} \ell^{\infty}_{m_0} \ell^{1}_{n^{\prime}} \ell^{1}_{m^{\prime}}.$ Thus by Proposition \ref{kernel bdd in bilinear case}, we get the desired result.
		
		\begin{rem}
			Let $\sigma$ be the symbol satisfying \eqref{symbol estimate for continuity result}, and $N_1,N_2,N_3$ are same as in Theorem \ref{bdd of symbol in modulatiobn space}. In the view of estimate \eqref{estimate of symbol in bilinear FIO case}, we observe that if we choose the phase function $\Phi$ in such a way that the infinite matrix of the correponding Fourier integral operator satisfying the condition $(iii)$ of Proposition \ref{kernel bdd in bilinear case}, then $T_\sigma$ defined by \eqref{definition of multilinear}, is a bounded linear operator from $\tilde{M}_{\mu}^{p, q} \otimes \tilde{M}_{\mu}^{p, q}$ to $\tilde{M}_{\mu v_{-m_{1}-m_{2}-2 N_{3},-m_3}}^{p, q}.$
		\end{rem}

	\end{proof}	
	
\end{document}